# Non-exchangeable evolutionary and mean field games and their applications

Hidekazu Yoshioka[1, *], Motoh Tsujimura[2], Tomohiro Tanaka[3]

[1] Graduate School of Advanced Science and Technology, Japan Advanced Institute of Science and Technology, 1-1 Asahidai, Nomi, Ishikawa, Japan

[2] Faculty of Commerce, Doshisha University, Karasuma-Higashi-iru, Imadegawa-dori, Kamigyo-ku, Kyoto 602-8580, Japan

[3] Disaster Prevention Research Institute, Kyoto University, Gokasho, Uji, Kyoto 611-0011, Japan

[*] Corresponding author: yoshih@jaist.ac.jp, ORCID: 0000-0002-5293-3246

**Abstract**

A replicator dynamic for non-exchangeable agents in a continuous action space is formulated and its well-posedness is proven in a space of probability measures. The non-exchangeability allows for the analysis of evolutionary games involving agents with distinct (and possibly infinitely many) types. We also explicitly connect this replicator dynamic to a stationary mean field game, which determines the pairwise actions of the heterogeneous agents. Moreover, as a byproduct of our theoretical results, we show that a class of nonlinear voter models, recently the subject of increasing interest, called $q$-voter models, can be viewed as a replicator dynamic driven by a utility that is a power of the probability density. This implies that non-exchangeable and/or mean-field game formulations of these models can also be constructed. We also present computational examples of evolutionary and mean field game models using a finite difference method, focusing on tragedy of the commons and the $q$-voter model with non-exchangeable agents, of which are interesting cases from theoretical and computational perspectives.

***Statements & Declarations***

**Fundings** This study was supported by the Japan Society for the Promotion of Science (KAKENHI No. 25K00240).

**Competing interests** The authors declare no relevant financial or non-financial interests.

**Data availability** The data will be made available upon reasonable request to the corresponding author.

**Acknowledgments** N.A.

**Declaration of generative AI in scientific writing** The authors did not use generative AI in the writing of this manuscript.

## 1. Introduction

### 1.1 Research background

Interactions among agents shape the social dynamics in human society, such as opinion formation [63], green travel behavior [50], immigration between countries [25], political polarization [80], and epidemic spread [70]. Social dynamics are also important in understanding the behavior of non-human animals, including beehive collective dynamics [73], bird migration [77], and animal foraging behavior [55].

Modeling and analyzing the link between agent-level interactions and decision-making—and the resulting collective dynamics—have proven effective for deepening our understanding of real-word social phenomena. Particularly, evolutionary game models that describe the time evolution of actions (or opinions) of agents (e.g., Sandholm [67]) have been a central mathematical tool for this purpose. Most evolutionary game models take the form of ordinary differential equations (ODEs), partial differential equations, or their time-discrete versions. These equations govern the distribution of actions, driven by some pairwise protocols based on optimizing a utility to be maximized or minimized. A pairwise protocol defines how each agent updates their decisions in response to others' actions. Key examples include the replicator dynamic (RD), where agents respond to utility differences [26, 72]; the logit dynamic, where decisions incorporate uncertainties [47, 21], and the projection dynamic, where non-smooth effects influence decision-making process of agents. Cheung [19] introduced pairwise comparison dynamics, a generalization of these protocols, and analyzed their well-posedness. These dynamics have since been further studied in the context of system stability [14, 36, 37]. Nonlocal protocols, where agents consider the actions of multiple others (i.e., higher-order comparisons), have also been explored [5, 9]. Evolutionary game models have found applications in machine learning, such as similarity-based learning [54], no-regret learning [49], and multi-agent optimization via gradient descent tailored for multi-agent optimization [11]. Advanced topics include stochastic versions of RDs [78] and evolutionary games with mutation [22]. Applications to transboundary water management have been investigated as well [32, 42, 56]. A protocol directly connecting the flow and use of a common resource has been proposed in López-Corona et al. [51] focusing on the groundwater management as an example.

In the real world, agents in a society or a group are not homogeneous but have some degree of heterogeneity, motivating the development of evolutionary game models involving heterogeneous (also called non-exchangeable in literature) agents. Examples include models with age structures [2, 34], multi-population opinion dynamics [7], and RD with group selection [45]. Evolutionary game models with discrete agent types have been studied in the context of finite-action cases [27, 38]. Modeling agent heterogeneity has also been a major topic in non-exchangeable mean field models, which describe the macroscopic evolution of quantities defined on a (hyper)graph [6, 23, 52]. However, evolutionary game models incorporating heterogeneity have been studied less compared to the homogeneous ones possibly because of the increased complexity in both theory and application. In particular, their connection to mean field games (MFGs), as explained below, has not been well addressed to the best of our knowledge.

The MFG is an optimal control model in which agents make decisions based on future predictions, which are typically represented as expectations conditioned on current information [e.g., 30, 48]. This contrasts with evolutionary game models, where decisions are solely made based on the current state. MFG models have been extensively studied in various contexts, including energy management of hybrid electric vehicles [18], disease spread and control [62], asset liquidation in financial markets [87], and incentive design [69]. MFGs with heterogeneous agents have been analyzed in Markov decision-type models [24], renewable energy supply-demand balancing [64], and two-player market design games for reducing greenhouse gas emission [79] and for solar energy pricing [74].

Recently, it has been shown that an MFG model can reduce to an evolutionary game model under certain limit conditions where the discount rate, a parameter that determines the agent's degree of myopia, becomes large. Bertucci et al. [13] studied this in a diffusive agent-based model, and Bardi and Cardaliaguet [8] analyzed an aggregation and flocking model. Later, MFG versions of projection and replicator-type dynamics were addressed using an inverse control approach [83, 84]. From this perspective, problems involving many, or even a continuum of, non-exchangeable agents have not been extensively studied. Moreover, studying other types of evolutionary game models would further deepen our understanding of the link between evolutionary and mean field games.

### 1.2 Aim and contribution

This study has two aims. The main aim is to mathematically and computationally analyze an evolutionary game model involving non-exchangeable agents and its associated MFG version. The secondary aim is to show that recently studied agent-based models, known as $q$-voter models, can be reformulated as evolutionary games, thereby enabling the derivation of their corresponding MFG versions. We focus on the RD, because this is one of the most employed pairwise comparison protocols in evolutionary games. Our main contributions are as follows.

We focus on the RD on a continuous action space [60] with agents' heterogeneity—referred to as type—parameterized continuously. Under this setting, the joint domain of action and type is a closed cube, i.e., a compact set. Compactness is a standard assumption in the literature [e.g., 5, 53]. Our dynamic is mathematically a version of the model studied in Mendoza-Palacios and Hernández-Lerma [53], though it has not yet been analyzed for specific problem settings. Cases involving discrete actions and/or types can also be handled within this framework by considering measure-valued solutions to the RD. We prove that this RD is well-defined and admits a unique solution in a space of probability measures equipped with the total variation (TV) norm. The key assumption for well-posedness is the boundedness and Lipschitz continuity of the utility with respect to the probability measures.

We are also interested in the applicability of the pairwise comparison protocol. We show that $q$-voter models [16] can be reformulated as an RD, with the utility given by a power function of the probability density. A $q$-voter model is a system of ODEs inspired by opinion dynamics with social interactions, with the "$q$" in its name indicating the degree of nonlinearity in these interactions. These models and related ones have been extensively studied, focusing on system stability and macroscopic ordering properties [e.g., 39, 59, 65, 66]. This study provides a new perspective on $q$-voter models from the standpoint of pairwise protocols because they have not been previously considered in a continuous action space. The $q$-voter version of the RD presents a challenging case where the Lipschitz continuity and boundedness of the corresponding utility are not satisfied. We propose regularization methods to address this issue.

We link the evolutionary game model and its MFG version through an inverse control approach, where a pairwise protocol—identified as the jump rate of agents' actions—is formally derived by choosing a specific objective function to be optimized using a dynamic programming principle. The optimized objective function, called the value function, corresponds to the utility in the evolutionary game. This methodology was employed in previous studies for time-dependent MFGs [83, 84]; however, the MFGs derived in these studies were investigated only computationally. In this paper, we study a stationary MFG in an infinite horizon (Chapter 7 in Bensoussan et al. [10]), which has been analyzed mathematically and computationally [1, 12] but not in the context of evolutionary games. With this formulation, we show, using a fixed-point argument, that our MFG admits a unique solution—a pair consisting of a probability measure and a value function, which possesses some regularity, sign, and boundedness properties.

Finally, we present computational applications of the proposed evolutionary and mean field games, which are discretized using a finite difference method. We consider several cases covered by the proposed mathematical models, including the aggregative game [20, 43, 44] and the $q$-voter model with a continuous action space. As computational examples, we first investigate public goods games related to the tragedy of the commons, where many agents exploit the same resource. While the problems are relatively simple, the non-exchangeability of agents leads to non-trivial grouping depending on the discount rate. We also analyze a utility that potentially prevents the presence of fully exploiting agents. The $q$-voter model with non-exchangeable agents is investigated for the first time in this paper.

The rest of this paper is organized as follows. **Section 2** presents our RD and examines its well-posedness. It also discusses the adaptation of the $q$-voter model to this dynamic. **Section 3** derives and presents the MFG version of the RD, and the well-posedness of the MFG is also studied. **Section 4** presents computational experiments of both the evolutionary and MFG models. **Section 5** concludes the paper and provides perspectives on this study. **Appendix** contains auxiliary results (**Section A.1**) and proofs of the propositions (**Section A.2**).

## 2. Replicator dynamic

We present and briefly analyze our RD.

### 2.1 Preliminaries and model

We introduce several notations. The explanations here are based on previous studies of evolutionary and dynamic games with continuous action spaces [e.g., 6, 47, 52, 60, 84].

Time $t \geq 0$ is a continuous parameter. The space of the agents' actions is given by the compact set $\Omega = [0,1]$. The space of the agent's types is given by the compact set $I = [0,1]$. The two sets $\Omega$ and $I$ are equivalent, but we distinguish them because they represent different objects. One may replace $\Omega$ or $I$ by a compact interval with proper scaling and shifts.

The Borel $\sigma$-field on the domain $\Omega$ is denoted as $\mathcal{B}$. The space of all finite signed measures in the measurable space $(\mathcal{B},\Omega)$ is denoted as $\mathcal{M}(\mathcal{B},\Omega)$. We equip the space $\mathcal{M}(\mathcal{B},\Omega)$ with the TV norm $\|\cdot\|_{\mathrm{TV}}$, making it a Banach space. Here, $\|\mu\|_{\mathrm{TV}}=\sup_{f}\left|\int_{\Omega}f(x)\mu(\mathrm{d}x)\right|$ for any $\mu\in\mathcal{M}(\mathcal{B},\Omega)$, where the supremum is taken over all bounded and measurable functions such that $\sup_{x\in\Omega}|f(x)|\le 1$. The collection of all signed measures with a TV norm not greater than 2 is denoted as $\mathcal{M}_2(\mathcal{B},\Omega)=\left\{\mu\in\mathcal{M}\,\middle|\,\|\mu\|_{\mathrm{TV}}\le 2\right\}$. The collection of all probability measures in $(\mathcal{B},\Omega)$ is denoted as $\mathcal{P}(\mathcal{B},\Omega)$. In the sequel, we omit the notation "$(\mathcal{B},\Omega)$" from the aforementioned function spaces for simplicity. When a measure $\mu\in\mathcal{M}$ is time-dependent, it is expressed as $\mu_t$ at time $t$. The collection of all continuous mappings from a space $Q_1$ to a space $Q_2$ is denoted as $C(Q_1,Q_2)$.

We also define a continuum of signed measures on $\Omega\times I$ as $\{\mu(y,\mathrm{d}x)\}_{y\in I}$, i.e., $\mu(y,\cdot)\in\mathcal{M}$ for any $y\in I$. The collection of all such signed measures is denoted as $\mathcal{M}^{(I)}$, to which we associate the norm $\|\mu\|_{\mathrm{TV}}^{(I)}=\sup_{y\in I}\|\mu(y,\cdot)\|_{\mathrm{TV}}$. The continuum versions of $\mathcal{M}_2$ and $\mathcal{P}$ are defined as $\mathcal{M}_2^{(I)}=\left\{\mu\in\mathcal{M}^{(I)}\,\middle|\,\|\mu\|_{\mathrm{TV}}^{(I)}\le 2\right\}$ and $\mathcal{P}^{(I)}=\left\{\mu\in\mathcal{M}^{(I)}\,\middle|\,\mu(y,\cdot)\in\mathcal{P},\ y\in I\right\}$, respectively. If $\mu\in\mathcal{M}^{(I)}$ depends on time $t$, we write it as $\{\mu_t(y,\mathrm{d}x)\}_{t\ge 0}$ or simply $\mu_t$ when the meaning is clear from the context.

In this paper, a utility $U$ is a mapping from $\Omega\times I\times\mathcal{M}^{(I)}\to\mathbb{R}$ that satisfies the following three conditions:

$$0\le U(x,y,\mu)\le\bar{U}\quad\text{for any }(x,y,\mu)\in\Omega\times I\times\mathcal{M}_2^{(I)},\tag{1}$$

$$|U(x_1,y_1,\mu)-U(x_2,y_2,\mu)|\le\phi(|x_1-x_2|+|y_1-y_2|)\quad\text{for any }(x_1,y_1,x_2,y_2,\mu)\in\Omega\times I\times\Omega\times I\times\mathcal{M}_2^{(I)},\tag{2}$$

$$|U(x,y,\mu)-U(x,y,\nu)|\le L_U\|\mu-\nu\|_{\mathrm{TV}}^{(I)}\quad\text{for any }(x,y,\mu,\nu)\in\Omega\times I\times\mathcal{M}_2^{(I)}\times\mathcal{M}_2^{(I)}.\tag{3}$$

Here, $\bar{U}>0$ is a constant that bounds the size of utilities, $L_U>0$ is the Lipschitz constant of utilities, and $\phi$ is a non-negative and non-decreasing continuous function such that $\phi(0)=0$. The conditions (1), (2), and (3) represent boundedness, uniform continuity for parameters, and Lipschitz continuity for measures. These conditions are central to well-defining the RD and its associated MFG.

Our RD is the following formal differential equation that governs a time-dependent probability measure $(\mu_t)_{t\ge 0}$ with the spatial domain $\Omega\times I$:

$$\underbrace{\frac{\mathrm{d}\mu_t(y,\mathrm{d}x)}{\mathrm{d}t}}_{\text{Temporal change}}=\left(\underbrace{U(x,y,\mu_t)}_{\text{Current utility}}-\underbrace{\int_{z\in\Omega}U(z,y,\mu_t)\mu_t(y,\mathrm{d}z)}_{\text{Averaged utility}}\right)\underbrace{\mu_t(y,\mathrm{d}x)}_{\text{Rate change}}\tag{4}$$

for all $(x,y,t) \in \Omega \times I \times (0,+\infty)$ subject to an initial condition $\mu_0 \in \mathcal{P}^{(I)}$. More rigorously, the RD (4) should be understood for any Borel measurable sets $A \in \mathcal{B}$ with respect to $\mathrm{d}x$ (i.e., "$\mathrm{d}x$" in (4) is replaced by "$A$"). According to (4), $\mu$ evolves over time by comparing the current and averaged utilities, and agents with distinct types interact through their utilities.

As pointed out in Cheung [19], a class of evolutionary game models including RDs, can be expressed as mass conservation laws. In our context, the RD (4) can alternatively be expressed as follows, as can be verified analytically:

$$\begin{aligned} \underbrace{\frac{\mathrm{d}\mu_t(y,\mathrm{d}x)}{\mathrm{d}t}}_{\text{Temporal change}} &= \underbrace{\int_{z\in\Omega}\left(U(x,y,\mu_t)-U(z,y,\mu_t)\right)_+ \mu_t(y,\mathrm{d}z)\mu_t(y,\mathrm{d}x)}_{\text{Probability inflow}} \\ &\quad -\underbrace{\int_{z\in\Omega}\left(U(z,y,\mu_t)-U(x,y,\mu_t)\right)_+ \mu_t(y,\mathrm{d}z)\mu_t(y,\mathrm{d}x)}_{\text{Probability outflow}}, \end{aligned} \tag{5}$$

where $(a)_+ = \max\{a,0\}$ for any $a \in \mathbb{R}$. The expression (5) implies that the utility difference drives the decision-making process of agents, where agents choose actions that yield a higher average utility by comparing their current action $x$ and the other actions $z$. An important finding from (5) is that the following coefficient $\rho(x,y,\mathrm{d}z)$ becomes the jump intensity measure from state $x$ to $z$:

$$\rho(x,y,\mathrm{d}z) = \left(U(z,y,\mu_t)-U(x,y,\mu_t)\right)_+ \mu_t(y,\mathrm{d}z). \tag{6}$$

The actions of agents are therefore probably concentrated on those more heavily used. By (6), the net jump intensity from $x$ is given by $\int_{z\in\Omega}\rho(x,y,\mathrm{d}z)$ for each $x,y$. Thus, the RD (5) can be rewritten as

$$\frac{\mathrm{d}\mu_t(y,\mathrm{d}x)}{\mathrm{d}t} = \int_{z\in\Omega}\rho(z,y,\mathrm{d}x)\mu_t(y,\mathrm{d}z) - \int_{z\in\Omega}\rho(x,y,\mathrm{d}z)\mu_t(y,\mathrm{d}x). \tag{7}$$

We interchangeably use the forms (4), (5), and (7) in the sequel.

In this setting, any solution $\mu_t$ to (5) is a probability measure for all $y \in I$ and $t \geq 0$ (see **Proposition 1**). This also implies that the total mass is one, as we formally have

$$\frac{\mathrm{d}}{\mathrm{d}t}\int_{y\in I}\int_{x\in\Omega}\mu_t(y,\mathrm{d}x)\mathrm{d}y = \int_{y\in I}\frac{\mathrm{d}}{\mathrm{d}t}\left(\int_{x\in\Omega}\mu_t(y,\mathrm{d}x)\mathrm{d}y\right) = 0 \tag{8}$$

and

$$\int_{y\in I}\int_{x\in\Omega}\mu_0(y,\mathrm{d}x)\mathrm{d}y = \int_{y\in I}1\mathrm{d}y = 1. \tag{9}$$

Additionally, if $\mu_t(y,\cdot)$ admits a probability density $p_t(y,\cdot)$, then (5) can be rewritten as follows (with some abuse of notation):

$$\begin{aligned} \frac{\mathrm{d}p_t(x,y)}{\mathrm{d}t} &= \int_{z\in D}\left(U(x,y,p_t)-U(z,y,p_t)\right)_+ p_t(z,y)\mathrm{d}z p_t(x,y) \\ &\quad -\int_{z\in D}\left(U(z,y,p_t)-U(x,y,p_t)\right)_+ p_t(z,y)\mathrm{d}z p_t(x,y). \end{aligned} \tag{10}$$

We conclude this subsection with the well-posedness result for the RD.

***Proposition 1***

*For each,* $T>0$*, the RD (4) with the initial condition* $\mu_0 \in \mathcal{P}^{(1)}$ *admits a unique solution* $\mu \in C\left(\left[0,T\right],\mathcal{P}^{(1)}\right)$.

**Proof of this Proposition 1** will be found in **Section A.2**.

### 2.2 *q*-voter model

We establish a connection between the $q$-voter model and the RD. For simplicity, we consider the case where $\mu$ does not depend on $y$. The $q$-voter model is a system of ODEs that governs the time-dependent vector $\left\{\left(x_{i,t}\right)_{t\geq 0}\right\}_{i=1,2,3,\ldots,N}$ for some $N\in\mathbb{N}$ with $N\geq 2$ and $q>0$ with $q\neq 1$ (Eq. (29) of Ramirez et al. [65]):

$$\frac{\mathrm{d}x_{i,t}}{\mathrm{d}t}=\left\{x_{i,t}\right\}^{q}-x_{i,t}\sum_{j=1}^{N}\left\{x_{j,t}\right\}^{q} \tag{11}$$

for $i=1,2,3,\ldots,N$ and $t>0$ subject to an initial condition $\left\{x_{i,0}\right\}_{i=1,2,3,\ldots,N}\in\Delta_N$, where $\Delta_N$ is the simplex $\Delta_N=\left\{\mathbf{x}=\left\{x_i\right\}_{i=1,2,3,\ldots,N}:\text{each } x_i\geq 0,\ \sum_{i=1}^{N}x_i=1\right\}$. The equilibria of the $q$-voter model (11) have been classified in Ramirez et al. [65]; the solution converges to an interior point (often the center) of $\Delta_N$ when $q\in\left(0,1\right)$, and it diverges toward a boundary point when $q>1$. The parameter $q$ thus determines the stability of the model.

Because the solution to (11) is confined in $\Delta_N$, we can interpret each $x_i$ as the probability of the agents' action. More specifically, if we understand $x_i$ as the probability of actions in cell $C_i=\left[\left(i-1\right)\Delta x, i\Delta x\right)$ ($i=1,2,3,\ldots,N-1$) and $C_N=\left[\left(N-1\right)\Delta x,1\right]$, with a cell size $\Delta x=N^{-1}$, then (11) can be seen as a finite difference approximation of the RD of the form (10) (without the dependence on $y$) by appropriately scaling time. This motivates the following formal replacement: $x_{i,t}\to p_t\left(\hat{x}_i\right)$ $\sum_{j=1}^{N}\left\{x_{j,t}\right\}^{q}\to\int_{\Omega}\left\{p_t\left(z\right)\right\}^{q}\mathrm{d}z$, where $\hat{x}_i=\left(i-1/2\right)\Delta x$ is the center of the cell $C_i$. Then, (11) becomes

$$\frac{\mathrm{d}p_t\left(\hat{x}_i\right)}{\mathrm{d}t}=\left(\left\{p_t\left(\hat{x}_i\right)\right\}^{q-1}-\int_{\Omega}\left\{p_t\left(z\right)\right\}^{q-1}p_t\left(z\right)\mathrm{d}z\right)p_t\left(\hat{x}_i\right) \tag{12}$$

whose space-continuous version is the RD (10) with the power-type utility $U\left(x,p\right)=\left\{p\left(x\right)\right\}^{q-1}$. This linkage between replicator dynamics and $q$-voter models is novel and has not been found in the literature.

Interestingly, the continuous-action $q$-voter model presents challenges for both $q>1$ and $q\in\left(0,1\right)$. If $q>1$, the utility $U\left(x,p\right)=\left\{p\left(x\right)\right\}^{q-1}$ is Lipschitz continuous only if $q=2$, and in general, is unbounded because $p$ is unbounded as well. If $q\in\left(0,1\right)$, $U$ is unbounded even if $p$ is bounded.

Therefore, to computationally investigate the continuous-action $q$-voter model, we propose regularizing the utility as follows

$$U(x,p)=\begin{cases}\min\left\{p(x),\varepsilon^{-1}\right\}^{q-1} & \text{if } q>1\\ \max\left\{p(x),\varepsilon\right\}^{q-1} & \text{if } q\in(0,1)\end{cases}, \tag{13}$$

where $\varepsilon>0$ is a small regularization parameter. Moreover, the existence of the probability density $p_t$ is necessary for all $t>0$, suggesting that the initial condition $\mu_0$ must admit a density because $\mu_t$ ($t>0$) is absolutely continuous with respect to $\mu_0$ in RDs (e.g., Theorem 13 in Mendoza-Palacios and Hernández-Lerma [53]).

***Remark 1*** The right-hand side of the $q$-voter model reduces to zero if $q=1$, which is why we do not discuss this case.

## 3. Mean field game

### 3.1 Agent behavior

Following Yoshioka [84], we interpret the RD (10) as a Fokker–Planck (FP) equation that governs the continuum of mean-field SDEs describing the behavior of representative agents

$$\mathrm{d}X_t^{(y)}=\mathrm{d}Z_t^{(y)},\quad t>0 \tag{14}$$

subject to an initial condition $X_0^{(y)}\in\Omega$ for each $y\in I$, where $\left(Z_t^{(y)}\right)_{t\geq 0}$ is a continuum of càdlàg (right continuous with left limits) compound point processes with mutually-independent Poisson random measures (a Poisson random measure is a random measure that counts a number of jumps in a unit time). Each $Z_t^{(y)}$ has a jump intensity measure $\hat{u}_t(y,\mathrm{d}z)=\left(U(z,y,\mu_t)-U\left(X_{t-}^{(y)},y,\mu_t\right)\right)_+\mu_t(y,\mathrm{d}z)$, and a post-jump position $z$, i.e., the position after a jump is chosen from $\Omega$ according to the probability measure proportional to $\hat{u}_t(z,y)$. We aim to uncover an underlying optimality principle that leads to the SDE (14), to deepen our understanding of the evolutionary game model. Thus, our strategy is to derive an SDE of the form (14) in which the jump intensity is given by a measure-valued process $(u_t)_{t\geq 0}$, such that $u_t^*$ (the maximizing $u$ in the control problem presented below) is identified as $\hat{u}_t$. Phenomenologically, the SDE (14) can be motivated from consideration that each agent would determine their decisions based on other agents' actions in a probabilistic way, so that macroscopic, mean-field agent behavior emerges in an equilibrium. Specifically, what is chosen at each jump is the post-jump position but not jump size in our formulation.

### 3.2 Objective and optimality equations

We consider a specific objective function to relate the RD and an MFG. Let $\mu_t \in \mathcal{P}^{(I)}$ be the probability measure over actions of a representative agent who updates their action according on (14), where the jump intensity measure of $Z_t^{(y)}$ is formally given by $u_t = w_t m(y, \mathrm{d}z)$ with some non-negative scalar process $(w_t)_{t\geq 0}$ that is the variable to be optimized in the MFG, and $m$ is a time-averaged probability measure given by

$$m(y, \mathrm{d}x) = \delta \int_0^{+\infty} e^{-\delta s} \mu_s(y, \mathrm{d}x)\,\mathrm{d}s\,, \tag{15}$$

where $\delta > 0$ is the discount rate. We assume an initial condition $X_0^{(y)} = x \in \Omega$ for each $y \in I$. The process $(w_t)_{t\geq 0}$ is assumed to be non-negative, predictable (i.e., in our context, $w_t$ can be modelled using the left limit $X_{t-}^{(y)}$) with respect to the natural filtration generated by $Z^{(y)}$ ($y \in I$), and satisfies $\mathbb{E}\left[\int_0^{+\infty} e^{-\delta s} w_s^2 \mathrm{d}s\right] < +\infty$ (this condition is satisfied if $w_s$ is strictly bounded, which is the case for our optimal control in view of (23) and **Proposition 4**). For each $y \in I$, we consider a continuum of objective functions

$$J(x, y; w) = \mathbb{E}\left[\int_0^{+\infty} e^{-\delta s} \left\{ \underbrace{\delta U\left(X_s^{(y)}, y, m\right)}_{\text{Utility}} - \underbrace{C(w_s)}_{\text{Cost}} \right\} \mathrm{d}s \,\middle|\, X_0^{(y)} = x \right] \tag{16}$$

for any $(x, y) \in \Omega \times I$, and $C(w_s)$ is a non-negative cost function associated with activating jumps (i.e., changing actions).

In this formulation, we assume that the representative agent changes their action with an intensity that depends on the time-averaged probability measure $m(y, \mathrm{d}z)$. The use of $m$ indicates that agents make decisions based on their expectations of the near future distribution of actions. The objective function includes the cumulative utility and cost during the game, discounted over time. The discount factor $e^{-\delta s}$ is necessary to ensure a bounded value of $J$, and the scaling of the utility by $\delta$ follows the approach used in previous studies. This allows the optimized objective function to approach the utility $U$ in the limit $\delta \to +\infty$ where the agent becomes infinitely myopic and makes decisions based solely on the current state, as in evolutionary games [e.g., 83, 84]. The decision-making is based on both the current state and future predictions (the conditional expectation), with the timescale given by $\delta^{-1}$.

We also define the value function as a mapping $\Phi : \Omega \times I \to \mathbb{R}$:

$$\Phi(x, y) = \sup_{w} J(x, y; w). \tag{17}$$

Considering the form of the value function $\Phi$, the dynamic programming argument similar to that in Chapter 7 of Bensoussan et al. [10] shows that the governing equation for $\Phi$, also known as the Hamilton–Jacobi–Bellman (HJB) equation, should be formally given by:

$$\delta \Phi(x, y) = \delta U(x, y, m) + \sup_{w(\cdot)\geq 0} \left( \int_\Omega w(z) \left( \Phi(z, y) - \Phi(x, y) \right) m(y, \mathrm{d}z) - C(w) \right) \tag{18}$$

for any $(x,y)\in\Omega\times I$, where the supremum is taken over non-negative measurable functions. The maximizer of the right-hand side of (18) (if it exists) is denoted by $w^*(x,z,y)$ because the maximization is carried out at each $(x,z,y)$. The HJB equation should be coupled with the FP equation, which now reads

$$\frac{\mathrm{d}\mu_t(y,\mathrm{d}x)}{\mathrm{d}t}=\int_{z\in D}w^*(z,x,y)\mu_t(y,\mathrm{d}z)m(y,\mathrm{d}x)-\int_{z\in D}w^*(x,z,y)m(y,\mathrm{d}z)\mu_t(y,\mathrm{d}x). \quad (19)$$

We further specify $C(w)$ (that may depend on $y$) and reformulate the optimality equations (18) and (19) using an inverse control approach.

***Remark 2*** We have $m\in\mathcal{P}^{(I)}$ because of

$$\begin{aligned}\int_{x\in\Omega}m(y,\mathrm{d}x)&=\int_{x\in\Omega}\int_0^{+\infty}\delta e^{-\delta s}\mu_s(y,\mathrm{d}x)\mathrm{d}s\\&=\int_0^{+\infty}\delta e^{-\delta s}\left(\int_{x\in\Omega}\mu_s(y,\mathrm{d}x)\right)\mathrm{d}s\\&=\int_0^{+\infty}\delta e^{-\delta s}\mathrm{d}s\\&=1\end{aligned}. \quad (20)$$

## 3.3 Mean field game model

### 3.3.1 Formulation

We specify the cost function $C$ to complete the formulation. We choose

$$C(w_s)=\int_\Omega\frac{1}{2}w_s^2m(y,\mathrm{d}z), \quad (21)$$

which represents the quadratic cost that is proportional to the square of the jump intensity averaged with respect to the probability measure $m$. Thus, choosing a larger jump intensity or a more frequently used action on average increases the cost. Substituting (21) into (18) yields

$$\begin{aligned}&\sup_{w(\cdot)\geq 0}\left(\int_\Omega w(z)\left(\Phi(z,y)-\Phi(x,y)\right)m(y,\mathrm{d}z)-\int_\Omega\frac{1}{2}\left(w(z)\right)^2m(y,\mathrm{d}z)\right)\\&=\frac{1}{2}\int_\Omega\left(\Phi(z,y)-\Phi(x,y)\right)_+^2m(y,\mathrm{d}z)\end{aligned} \quad (22)$$

with the optimal $w=w^*$ as a maximizer being given by

$$w^*(x,z,y)=\left(\Phi(z,y)-\Phi(x,y)\right)_+. \quad (23)$$

Substituting $w=w^*$ into (19) yields the analogue of (5), where the utility $U$ is replaced by the value function $\Phi$. This similarity demonstrates that the present MFG is qualitatively linked to the RD, where the difference in value functions, rather than utility differences, drives the decision-making process of the representative agent.

Consequently, our MFG model is derived as follows: for any $(x,y)\in\Omega\times I$,

$$m(y,\mathrm{d}x)=\frac{\delta}{\delta-\left(\Phi(x,y)-\int_{z\in\Omega}\Phi(z,y)m(y,\mathrm{d}z)\right)}\mu_0(y,\mathrm{d}x)\ \left(=\mathbb{F}_\Phi(m)\right) \tag{24}$$

and

$$\Phi(x,y)=U(x,y,m)+\frac{1}{2\delta}\int_\Omega\left(\Phi(z,y)-\Phi(x,y)\right)_+^2 m(y,\mathrm{d}z)\ \left(=\mathbb{H}_m(\Phi)\right), \tag{25}$$

which will be further studied in the next subsection. Here, the (time-averaged version of) FP equation (24) is derived by multiplying $e^{-\delta t}$ against the both sides of (19) and integrating them over $(0,+\infty)$ for $t$.

***Remark 3*** The cost function (21) is a with-measure version of the one presented in Yoshioka [84], which assumed the existence of the density of $\mu$. This paper, however, demonstrates that assuming the existence of density is unnecessary.

***Remark 4*** Interestingly, our model is closely related to dynamical systems on graphons [3, 52]. In fact, the second term on the right-hand side of the HJB equation (25) encodes the connectivity between action pairs through $m$, which can be viewed as an integrable graphon. This, our MFG model can be interpreted as a type of graphon game. Note that this specific form of graphon games has not been studied before.

### 3.3.2 Mathematical analysis

We demonstrate that the system consisting of (24) and (25) admits a unique solution $(m,\Phi)$ in an appropriate setting. To achieve this, we define new two Banach spaces for studying the MFG: the space of bounded and continuous functions $\mathbb{X}=C(\Omega\times I)$ equipped with the usual $L^\infty$ norm $\|\cdot\|_{\mathbb{X}}$ on $\Omega\times I$, and the space of all continua of finite signed measures $\{m(y,\mathrm{d}x)\}_{y\in I}$, denoted as $\mathbb{Y}=C(I,\mathcal{M})$, equipped with the norm $\|\cdot\|_{\mathbb{Y}}=\|\cdot\|_{\mathrm{TV}}^{(I)}$, which forms a complete metric space because the space of continuous functions from one complete metric space to another is also complete [e.g., 28]. The spaces $\mathbb{X}$ and $\mathbb{Y}$ serve as solutions spaces of the HJB and FP equations, respectively. We set a space of probability measures $\mathbb{W}=C(I,\mathcal{P})\subset\mathbb{Y}$. We also set $\mathbb{X}_{\bar{U}}=\left\{\Phi\in\mathbb{X}:0\leq\Phi\leq\bar{U}\text{ on }\Omega\times I\right\}\subset\mathbb{X}$.

***Remark 5*** The mathematical analysis presented in this subsection is novel, even in case where there is no $y$-dependence.

We first present two results that deal with the FP and HJB equations separately.

***Proposition 2***

*Let* $\delta > 2\bar{U}$ *, and let* $l \in \mathbb{W}$ *be given. Any solution* $\Phi \in \mathbb{X}$ *to* $\Phi = \mathbb{H}_l(\Phi)$ *satisfies* $0 \le \Phi \le \bar{U}$ *on* $\Omega \times I$ *. Moreover, the equation* $\Phi = \mathbb{H}_l(\Phi)$ *admits a unique solution* $\Phi \in \mathbb{X}_{\bar{U}}$ *. The following estimate also holds:*

$$\left\| \Phi_1 - \Phi_2 \right\|_{\mathbb{X}} \le \left(1 - \frac{2\bar{U}}{\delta}\right)^{-1} \left(L_U + \frac{\bar{U}^2}{2\delta}\right) \left\| l_1 - l_2 \right\| = O(1) \left\| l_1 - l_2 \right\| , \tag{26}$$

*where* $\Phi_i$ *represents the solution to* $\Phi = \mathbb{H}_{l_i}(\Phi)$ *(* $i = 1,2$ *) and the right-most equality holds true for a sufficiently large* $\delta$ *.*

**Proof of this Proposition 2** will be found in **Section A.2**.

***Proposition 3***

*Assum* $\delta > \frac{3+\sqrt{5}}{2}\bar{U}$ *. Let* $\theta \in \mathbb{X}_{\bar{U}}$ *be given. Then, the equation* $m = \mathbb{F}_\Phi(m)$ *admits a unique solution* $m \in \mathbb{W}$ *. Moreover, the following estimate also holds:*

$$\left\| m_1 - m_2 \right\|_{\mathbb{Y}} \le \left(1 - \frac{\delta\bar{U}}{\left(\delta - \bar{U}\right)^2}\right)^{-1} \frac{2\delta}{\left(\delta - \bar{U}\right)^2} \left\| \theta_1 - \theta_2 \right\|_{\mathbb{X}} = O\left(\delta^{-1}\right) \left\| \theta_1 - \theta_2 \right\|_{\mathbb{X}} . \tag{27}$$

*Here, the unique solution to the equation* $m = \mathbb{F}_{\theta_i}(m)$ *with* $\theta = \theta_i \in \mathbb{X}_{\bar{U}}$ *is denoted as* $m_i$ *(* $i = 1,2$ *), and the right-most equality of (27) holds true for a sufficiently large* $\delta$ *.*

**Proof of this Proposition 3** will be found in **Section A.2**.

By **Propositions 2 and 3**, we conclude the following unique existence result for the MFG model (24) and (25).

***Proposition 4***

*Assume* $\delta > \frac{3+\sqrt{5}}{2}\bar{U}$ *. Then, the MFG model (24) and (25) admits a unique solution* $(m,\Phi) \in \mathbb{W} \times \mathbb{X}_{\bar{U}}$ *. Moreover, this* $m$ *is absolutely continuous with respect to* $\mu_0$ *.*

**Proof of this Proposition 4** will be found in **Section A.2**.

Finally, **Proposition 5** below demonstrates a consistency between the evolutionary game and the MFG in the limit as $\delta \to +\infty$ . This result implies that the two models coincide in this limit, suggesting that their differences are more pronounced when $\delta$ is smaller.

***Proposition 5***

*Let $\Phi$ and $m$ be the solution to the MFG model (24)-(25). Then, the limit $\lim_{\delta\to+\infty}(\Phi,m)=(\Phi_\infty,m_\infty)$ is a pair of a function $\Phi_\infty$ and a probability measure $m_\infty\in\mathbb{W}$ such that*

$$\Phi_\infty=U(\cdot,\cdot,\mu_0)\in\mathbb{X}_{\bar{U}} \text{ and } m_\infty=\mu_0\in\mathbb{W}. \tag{28}$$

**Proof of this Proposition 5** will be found in **Section A.2**.

***Remark 6*** The uniqueness result of our MFG model does not rely on any monotonicity condition on the utility $U$ or the operators involving it, which is a condition that has commonly been used to establish uniqueness in standard MFG theory [31, 40]. Instead, our analysis requires a sufficiently large discount $\delta$, which plays a role in regularizing the solutions.

# 4. Numerical computation

## 4.1 Discretization method

We discretize the MFG model (24) and (25) using a finite difference method similar to those in the literature [e.g., 4, 46, 84], combined with an iterative update method. We set the grid points $\mathrm{P}_{i,j}=(\hat{x}_i,\hat{y}_j)$ with $\hat{x}_i=(i-1/2)\Delta x$ ($i=1,2,3,...,N_x$) and $\hat{y}_j=(j-1/2)\Delta y$ ($j=1,2,3,...,N_y$) at which the two equations are discretized. Here, $N_x,N_y\in\mathbb{N}$ determine the computational resolutions, and we set $\Delta x=1/N_x$ and $\Delta y=1/N_y$. The discretized $m$ and $\Phi$ at $\mathrm{P}_{i,j}$ are denoted by $m_{i,j}$ and $\Phi_{i,j}$, respectively. For problems without $y$ dependence, it suffices to set $N_y=1$.

Equations (24) and (25) are discretized at each $\mathrm{P}_{i,j}$ as follows:

$$0=\delta\mu_{0,i,j}-\left(\delta-\left(\Phi_{i,j}-\sum_{k=1}^{N_x}\Phi_{k,j}m_{k,j}\right)\right)m_{i,j}\ \left(=\mathbb{F}_{i,j}(m,\Phi)\right) \tag{29}$$

and

$$0=-\Phi_{i,j}+U_{i,j}+\frac{1}{2\delta}\sum_{k=1}^{N_x}\left(\Phi_{k,j}-\Phi_{i,j}\right)_+^2 m_{k,j}\ \left(=\mathbb{H}_{i,j}(m,\Phi)\right), \tag{30}$$

where $U_{i,j}$ denotes the discretized utility $U$ at $\mathrm{P}_{i,j}$, and $\mu_{0,i,j}$ is the discretized $\mu_0$ at $\mathrm{P}_{i,j}$ with $\sum_{i=1}^{N_x}\mu_{0,i,j}=1$ ($j=1,2,3,...,N_y$). The full discretization method for $U_i$ may vary depending on the application and will be specified later. The probability density $p_{i,j}$ at each $\mathrm{P}_{i,j}$ can be recovered via $p_{i,j}=N_x m_{i,j}$.

We solve the coupled nonlinear equations (29) and (30) using the following naïve iteration method that resembles the forward-forward MFG employed in some applications [29, 75, 85]. In this method, the superscript $(n)$ ($n=0,1,2,...$) represents the numerical value at the $n$ th iteration. Here, $\Delta t>0$ is an increment of pseudo time, and $\varepsilon>0$ is a small error threshold value. We assume $\varepsilon<<\Delta t$.

| **Algorithm 1: Numerical computation of MFG** |
|---|
| 1: Step 1 Set an initial guess $\mu^{(0)}$ such that $\sum_{i=1}^{N_x} \mu_{i,j}^{(0)} = 1$ ( $j = 1,2,3,...,N_y$). Set an initial guess $\Phi^{(0)}$ such that all $\Phi_{i,j}^{(0)}$ are bounded between 0 and $\bar{U}$ . Set $n = 0$.<br>2: Step 2 Compute $\Phi_{i,j}^{(n+1)} = \Phi_{i,j}^{(n)} + \delta \Delta t \mathbb{H}_{i,j}\left(m^{(n)}, \Phi^{(n)}\right)$ at all grid points.<br>3: Step 3 Compute $m_{i,j}^{(n+1)} = m_{i,j}^{(n)} + \Delta t \mathbb{F}_{i,j}\left(m^{(n)}, \Phi^{(n)}\right)$ at all grid points.<br>4: Step 4 Compute the error $\max\left\{N_x \left\| m_{i,j}^{(n+1)} - m_{i,j}^{(n)} \right\|, \left\| \Phi_{i,j}^{(n+1)} - \Phi_{i,j}^{(n)} \right\|\right\}$ among all grid points, and their maximum value is denoted by $\varepsilon^{(n)}$.<br>5. Step 5 If $\varepsilon^{(n)} \le \varepsilon$ , then the couple $\left(m^{(n+1)}, \Phi^{(n+1)}\right)$ is chosen as the numerical solution and terminate the algorithm. If it is not, then set $n \to n+1$ and go to Step 2. |

We have the proposition stating that numerical solutions obtained from **Algorithm 1** are reasonable in the sense that they are non-negative, bounded, and conserve the probability as in our MFG model.

***Proposition 6***

*Assume $\delta > 2\bar{U}$ and set $K = \delta + \sqrt{\delta\left(\delta - 2\bar{U}\right)}$. Then, for each $n \in \mathbb{N}$, the following statements hold true for $\Delta t < \dfrac{1}{\delta + K}$ at all $n = 0,1,2,...$:*

$$m_{i,j}^{(n)} \ge 0 \quad \text{at all grid points (non-negativity),} \tag{31}$$

$$0 \le \Phi_{i,j}^{(n)} \le K \quad \text{at all grid points (boundedness),} \tag{32}$$

$$\sum_{i=1}^{N_x} m_{i,\cdot}^{(n)} = 1 \quad \text{(conservation of probability).} \tag{33}$$

**Proof of this Proposition 6** will be found in **Section A.2**.

We conclude the theoretical numerical analysis with a proposition stating that the system (29) and (30) is uniquely solvable. According to **Propositions 6** and **7**, we have that **Algorithm 1** potentially yields the unique numerical solution to the discretized system (29) and (30). Here, we write $\mathbf{m} = \left\{m_{i,j}\right\}_{\substack{1 \le i \le N_x \\ 1 \le j \le N_y}}$ .

***Proposition 7***

*Fix arbitrary* $N_x, N_y \in \mathbb{N}$. *Assume that the following inequality holds true for any* $\mathbf{m}^{(k)}$ *such that* $m_{\cdot,\cdot}^{(k)} \geq 0$, $\sum_{i=1}^{N_x} m_{i,\cdot}^{(k)} = 1$ *(*$k = 1,2$*): for all* $1 \leq i \leq N_x$ *and* $1 \leq j \leq N_y$*,*

$$\left| U_{i,j,1} - U_{i,j,2} \right| \leq L_{\bar{U}} \max_{1 \leq j' \leq N_y} \sum_{i'=1}^{N_x} \left| m_{i',j'}^{(1)} - m_{i',j'}^{(2)} \right|. \tag{34}$$

*Then, for a sufficiently large* $\delta > 0$*, the system (29) and (30) admits a unique solution. Here,* $U_{i,j,k} = U_{i,j}\big|_{\mathbf{m}=\mathbf{m}^{(k)}}$ *(*$k = 1,2$*). This solution consists of a pair of sequences* $\left\{ m_{i,j} \right\}_{\substack{1 \leq i \leq N_x \\ 1 \leq j \leq N_y}}$ *and* $\left\{ \Phi_{i,j} \right\}_{\substack{1 \leq i \leq N_x \\ 1 \leq j \leq N_y}}$ *such that* $0 \leq \Phi_{i,j} \leq \bar{U}$ *and* $m_{i,j} \geq 0$ *at all grid points and* $\sum_{i=1}^{N_x} m_{i,j} = 1$ *(*$1 \leq j \leq N_y$*).*

**Proof of this Proposition 7** will be found in **Section A.2**.

Finally, we consider that, at least technically, the discretization method presented in this paper is extendable to cases where the domain $\Omega$ is multidimensional.

## 4.2 The tragedy of the commons

### 4.2.1 Utility

As an application of the evolutionary and MFG models, we focus on the tragedy of the commons [33]. According to Ostrom [61], "*each herder received large benefits from selling his or her own animals while facing only small costs of overgrazing.*". This concept has been studied in various MFG-related contexts, including agent-based models [57], traveling waves in spatially distributed systems [41], Q-learning protocols [82], and models incorporating resource-agent feedback [17]. From a resource and environmental management perspective, related applications include pollution control [81], groundwater management [58], and sustainable tourism management [71]. However, MFGs involving non-exchangeable agents have not yet been explored in this context. Note that MFGs of the tragedy of commons case with exchangeable agents can be studied by simply reducing $I$ to a singleton, although it is out of interest in this paper because our focus is on the non-exchangeability that has not been studied well. In some application, a non-exchangeable case reduces to an exchangeable one in some case (e.g., see, **Fig. 7** presented later). We note here that being an exchangeable case may not mean that it is simpler than a non-exchangeable case. It will depend on utilities.

In our setting, the simplest form of the tragedy of the commons can be represented by the following utility, for each $(x, y, \mu) \in \Omega \times I \times \mathcal{P}^{(I)}$, we set

$$U(x, y, \mu) = \left( \underbrace{f(y) A(\alpha)}_{\text{Marginal benefit}} - \underset{\text{Unit cost}}{c} \right) \underset{\text{Harvesting}}{x} \quad \text{with } \alpha = \int_{\Omega \times I} x \mu(y, \mathrm{d}x) \mathrm{d}y. \tag{35}$$

This utility contains the benefit $f(y) A(\alpha) x$ and cost $cx$ for harvesting a common resource, where $f$ is a continuous function assumed to be increasing, and $A$ is a decreasing function of the average harvesting $\alpha$. The increasing nature of $f$ implies that an agent with a larger $y$ places more weight on the benefit than the cost. The decreasing nature of $A$ implies that an increase in the harvesting $x$ on average decreases the benefit or equivalently increases the relative cost. This model can be seen as a continuous-type analog of the abstract model in Cheung and Lahkar [20] for a game with continuous actions.

Our computation of the RD and MFG uses the common computational resolution of $N_x = N_y = 200$, the pseudo time increment of $\Delta t = 0.2$ except for the cases with $\delta = 10$ where we set $\Delta t = 0.02$ for a stability reason. The error threshold of $\varepsilon = 10^{-10}$. In our setting, **Algorithm 1** is terminated at most $O(10^4)$ iterations. The utility (35) at $\mathrm{P}_{i,j}$ is discretized as follows:

$$U_{i,j} = \left(f\left(\hat{y}_j\right)A\left(\hat{\alpha}\right) - c\right)\hat{x}_i \quad \text{with} \quad \hat{\alpha} = \sum_{j=1}^{N_y}\sum_{i=1}^{N_x}\hat{x}_i m_{i,j} . \tag{36}$$

The initial condition $\mu_0$ is the uniform distribution in $\Omega \times I$. In the sequel, by a (discounted version of the) RD, we mean

$$m(y,\mathrm{d}x) = \frac{\delta}{\delta - \left(U(x,y,m) - \int_{z\in\Omega} U(z,y,m)\, m(y,\mathrm{d}z)\right)} \mu_0(y,\mathrm{d}x) . \tag{37}$$

Its well-posedness follows in the same way as **Proposition 3**. This equation is numerically discretized following the procedure in **Section 4.1** and is solved by **Algorithm 1**, where $\Phi$ is replaced by $U$.

### 4.2.2 Computational results

First, we examine the test case with $f(y) = 1 + 2y$, $A(\alpha) = \alpha^{-2}$, and $c = 4$. This case does not satisfy the Lipschitz continuity condition (3), but it can be regularized without affecting the numerical solutions by replacing $A(\alpha)$ with $\max(\alpha, \Delta x)^{-2}$ (one may use a sufficiently small constant instead of $\Delta x$), which we use in the sequel. The coefficients and parameter values in this computational case are chosen such that the solution profiles visibly change across the values of $\delta$.

In light of **Proposition 5**, we compare the probability densities $p$ between the RD and MFG, as shown in **Fig. 1** and **Table 1** for different values of $\delta$. We also compare the utility $U$ of the RD and value function $\Phi$ of the MFG as shown in **Fig. 2** and **Table 2** for different values of $\delta$. Here, we quantify the difference between the two models by the maximum values among and average value at all the grid points. The computational results suggest that the difference in the density $p$ is of the order $\delta^{-1}$, quantifying the convergence rate that is not explicitly shown in **Proposition 5**. For the utility $U$ of the RD and $\Phi$ of the MFG, the convergence speed between them is slower than that of the density $p$, but it still decreases with respect to $\delta$.

Profiles of the numerical solutions, particularly the density $p$, differ qualitatively between large and small values of $\delta$, as shown in **Fig. 1**. For large $\delta$, the profile of $p$ is bimodal, with concentrations around $(x, y) = (0, 0)$ and $(x, y) = (1, 1)$, suggesting that agents who place greater weight on the benefits of harvesting tend to exploit the resource, while who place less weight refrain from doing so. This behavior can be interpreted as the emergence of the tragedy of the commons. These characteristics become more pronounced for smaller $\delta$, where agents are almost fully separated into two distinct groups: one that fully exploits the resource ($x \approx 1$ and large $y$) and another that fully refrains from it ($x \approx 0$ and small $y$). These computational results suggest that a continuum of non-exchangeable agents can, in effect, be represented by a finite number of representative types when decision-making is based on a long-term perspective ($\delta$ is small). This observation implies that some additional mechanism should be incorporated into the utility (35), as discussed in the next subsection.

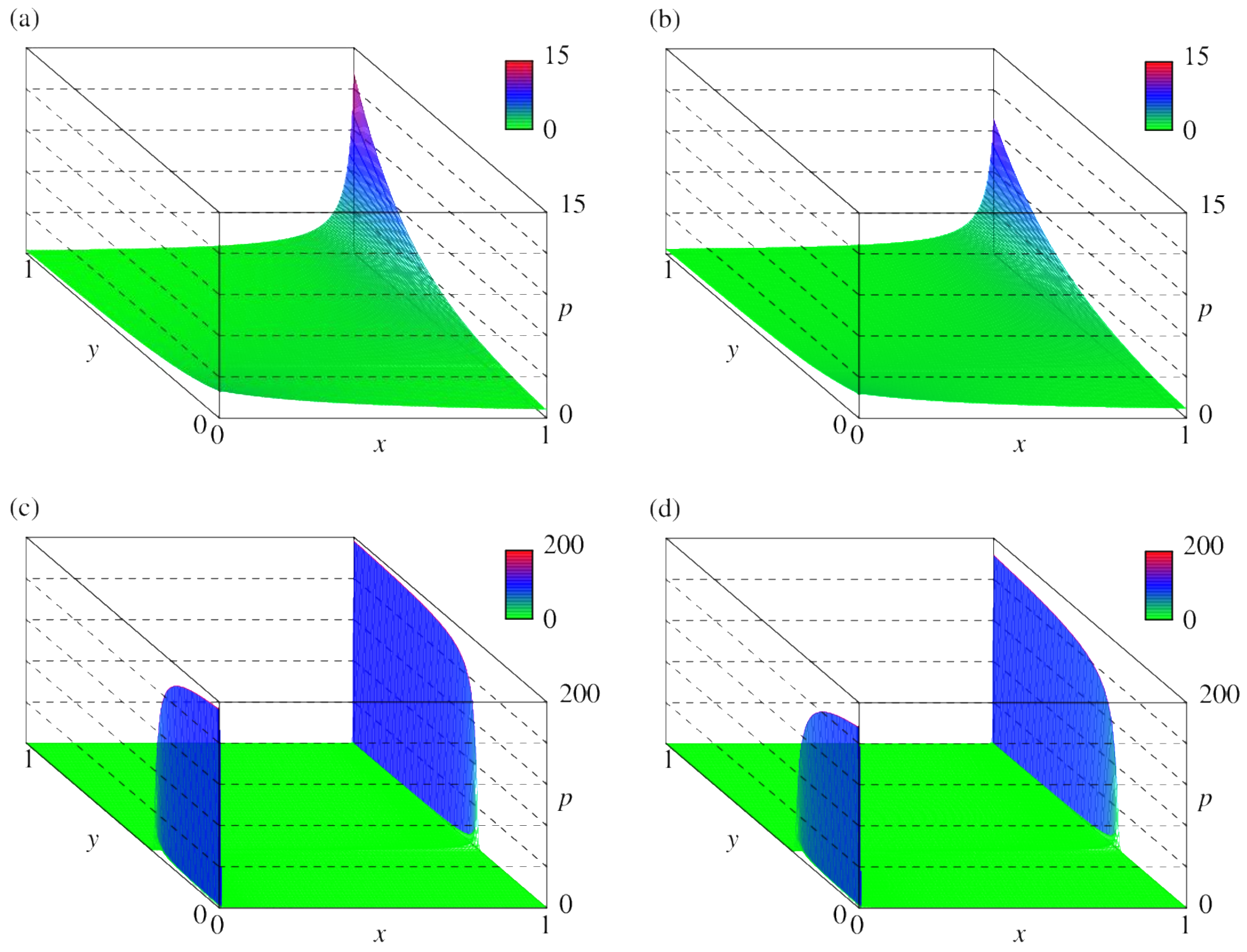

**Fig. 1** Computed probability density $p$. (a) RD with $\delta = 1$, (b) MFG with $\delta = 1$, (c) RD with $\delta = 0.01$, and (d) MFG with $\delta = 0.01$.

**Table 1** Maximum and average absolute differences in the probability density $p$ between the RD and MFG.

| $\delta$ | Maximum difference | Average difference |
| --- | --- | --- |
| 0.01 | $5.826 \times 10^{1}$ | $2.125 \times 10^{-1}$ |
| 0.1 | $3.597 \times 10^{1}$ | $2.692 \times 10^{-1}$ |
| 1 | $3.316 \times 10^{0}$ | $6.433 \times 10^{-2}$ |
| 10 | $2.880 \times 10^{-2}$ | $4.237 \times 10^{-3}$ |

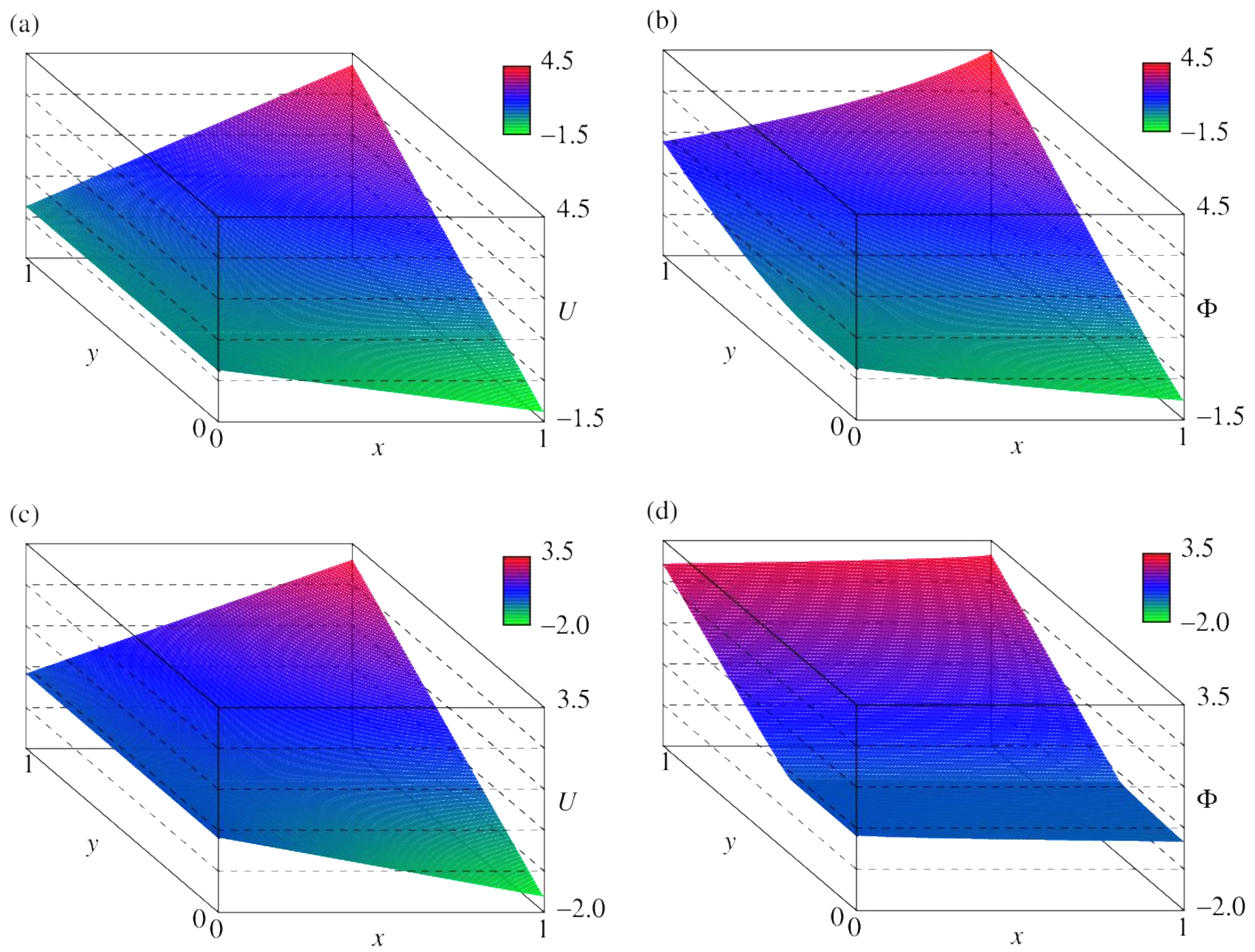

**Fig. 2** Computed utility $U$ or value function $\Phi$. (a) $U$ of RD with $\delta = 1$, (b) $\Phi$ of MFG with $\delta = 1$, (c) $U$ of RD with $\delta = 0.01$, and (d) $\Phi$ of MFG with $\delta = 0.01$.

**Table 2** Maximum and average absolute differences between the utility $U$ of RD and the value function $\Phi$ of MFG.

| $\delta$ | Maximum difference | Average difference |
|---|---|---|
| 0.01 | $2.864 \times 10^{0}$ | $5.743 \times 10^{-1}$ |
| 0.1 | $2.511 \times 10^{0}$ | $4.583 \times 10^{-1}$ |
| 1 | $1.814 \times 10^{0}$ | $2.958 \times 10^{-1}$ |
| 10 | $6.866 \times 10^{-1}$ | $7.967 \times 10^{-2}$ |

### 4.3 Tragedy of the commons with penalty

In the previous computational case, we found that the agents with higher $y$ values tend to choose actions close to full exploitation $x = 1$, triggering the tragedy of the commons. To address this, we consider an extended case in which a penalty is introduced into the utility:

$$U(x, y, \mu) = \Big( \underbrace{f(y) A(\alpha)}_{\text{Marginal benefit}} - \underbrace{c}_{\text{Unit cost}} \Big) \underbrace{x}_{\text{Harvesting}} - \underbrace{\begin{cases} P & \text{if } x \geq \underline{x} \\ 0 & \text{if } x < \underline{x} \end{cases}}_{\text{Penalty}} \quad \text{with } \alpha = \int_{\Omega \times I} x \mu(y, \mathrm{d}x) \mathrm{d}y . \tag{38}$$

Here, $P > 0$ denotes the strength of penalty, and $\underline{x} \in (0,1)$ is a threshold value such that choosing an action no smaller than it activates the penalty. A similar kind of discontinuous penalization has been used in Gao et al. [27]. In our context, the challenge with the utility (38) is that it is a discontinuous function of $x \in \Omega$, and hence does not satisfy the continuity condition (2), which was fundamental for the mathematical analysis in **Section 3**. We computationally investigate the RD and MFG under the discontinuous utility (38) for different values of $\delta$ and $P$. We fix $\underline{x} = 3/4$ here but choosing other values of $\underline{x}$ around this does not significantly change the discussion below. The discretization of utility (38) is conducted in the same manner as (36), where $x$ in the last term of (38) is replaced by $\hat{x}_i$ at each grid point. We did not encounter wiggles of numerical solutions with this methodology.

**Figs. 3**, **4**, and **5** compare the probability densities $p$ of RD and MFG for penalty values $P$ of 0.2, 1, and 5, respectively. Increasing the penalty strength $P$ more strongly suppresses the presence of fully exploiting agents. A small penalty does not significantly alter the probability density $p$ compared to the no-penalty case discussed in the previous subsection. When the penalty $P$ is sufficiently large, the population concentrates around the threshold value $\underline{x}$, likely because harvesting at a rate higher than $\underline{x}$ leads to a sudden increase in the net cost. Despite the discontinuity in the utility, numerical results show convergence of the MFG to the corresponding RD, as demonstrated in the **Appendix**. The computational results obtained suggest that the regularity of the utility $U$, particularly its defect, critically affects the agents' actions. From a computational standpoint, we can choose a more complex $U$ as long as $\delta$ is sufficiently large (**Proposition 7**); however, the well-posedness of the MFGs in such a case is still open.

We investigated how to resolve the tragedy of the commons by designing the utility, although this goal can also be achieved by modifying the dynamics themselves. For example, in socio-hydrology, Yu et al. [86] discussed how flood resilience arises from collective actions driven by cooperation among agents with memories of past flood events. This approach has since been applied to spatially explicit models of human-flood interaction [76]. From this perspective, enhancing the dynamics within the proposed mathematical framework would increase its applicability to real-world problems. A feasible step in this direction would be coupling the agent dynamics with an environmental dynamic, which is more straightforward in the RD framework, but more complex in the MFG because of the forward-backward nature of the equations. Incorporating such system dynamics into the RD or MFG could allow the model to determine the threshold $\underline{x}$ intrinsically, based on an equilibrium condition that balances the cost and benefit of harvesting. The proposed framework may thus serve as a foundation for addressing these more advanced research questions.

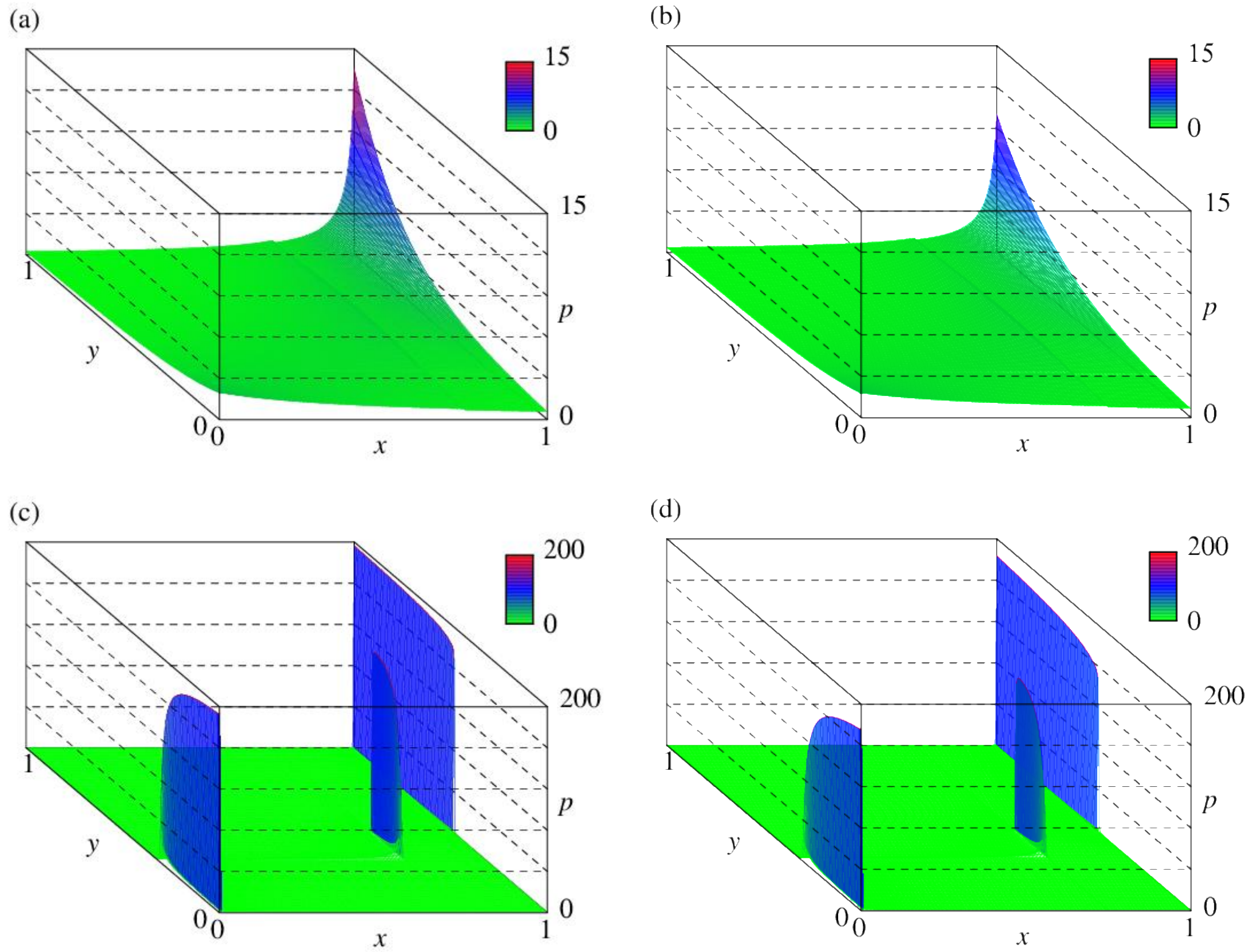

**Fig. 3** Computed probability density $p$ with $P = 0.2$. (a) RD with $\delta = 1$, (b) MFG with $\delta = 1$, (c) RD with $\delta = 0.01$, and (d) MFG with $\delta = 0.01$.

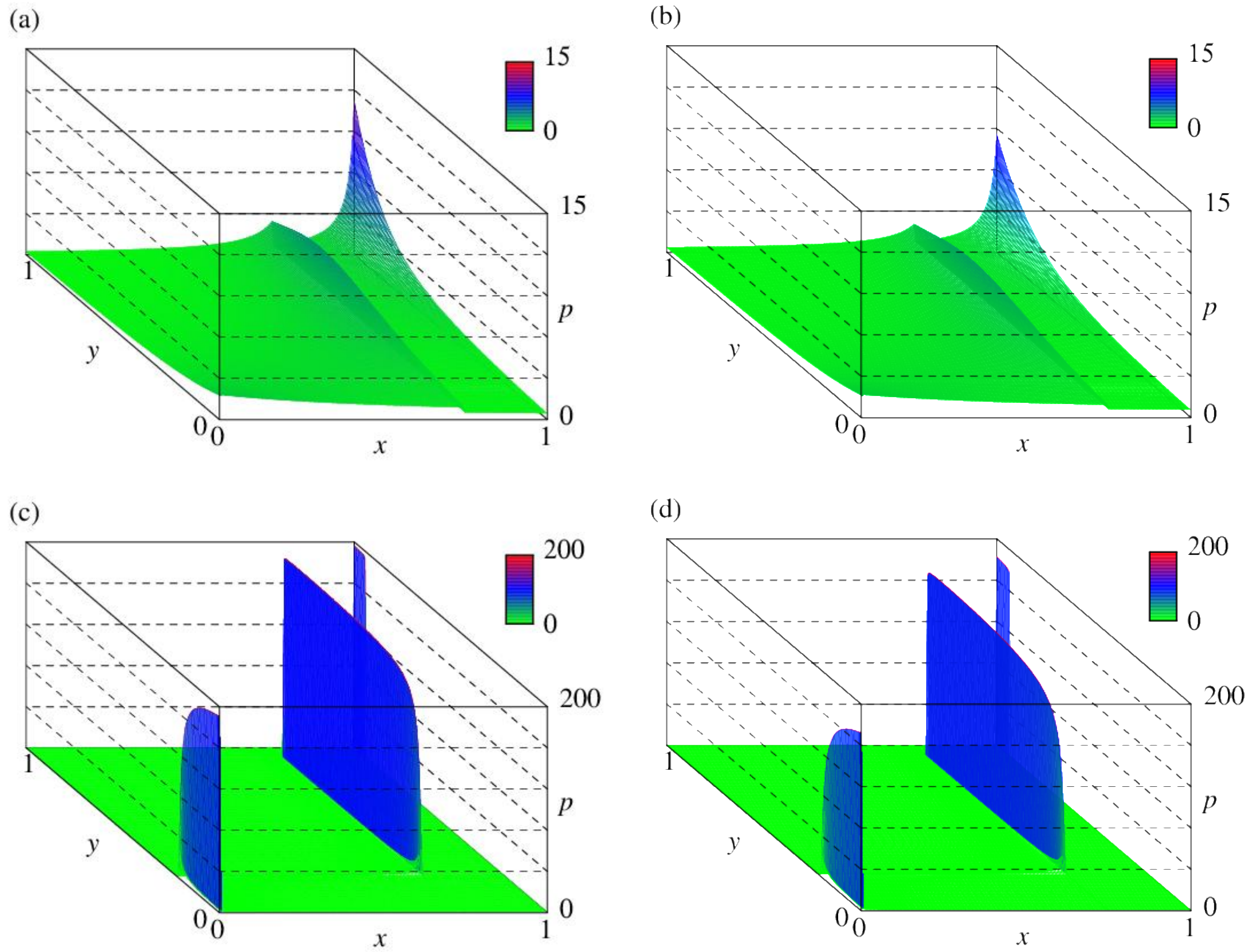

**Fig. 4** Computed probability density $p$ with $P = 1$. (a) RD with $\delta = 1$, (b) MFG with $\delta = 1$, (c) RD with $\delta = 0.01$, and (d) MFG with $\delta = 0.01$.

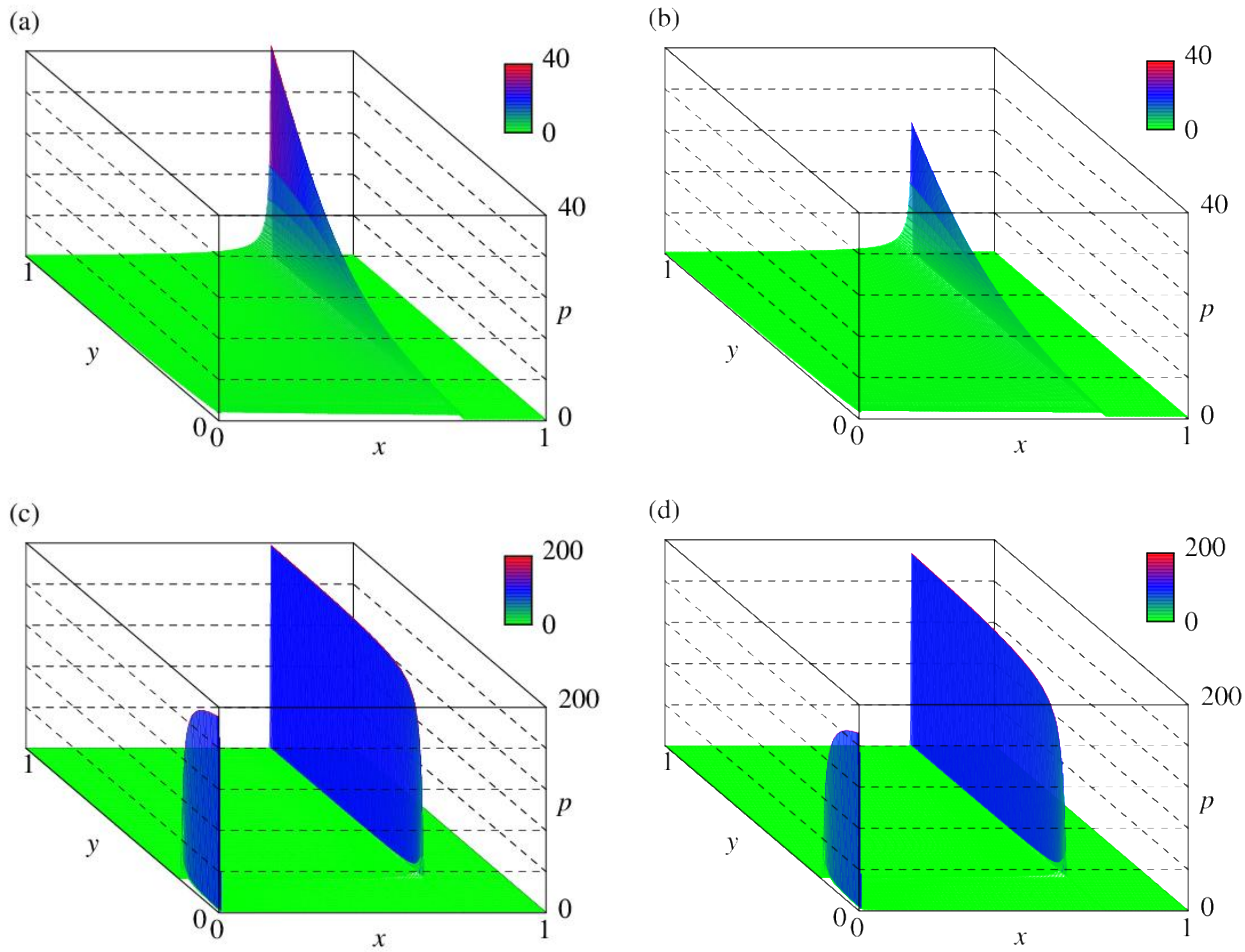

**Fig. 5** Computed probability density $p$ with $P = 5$. (a) RD with $\delta = 1$, (b) MFG with $\delta = 1$, (c) RD with $\delta = 0.01$, and (d) MFG with $\delta = 0.01$.

### 4.4 *q*-voter model

We also computationally apply the *q*-voter model to our MFG framework because it presents an interesting case where the utility is irregular and the assumptions required for each proposition are not satisfied. We consider both exchangeable and non-exchangeable cases.

For the exchangeable case, we can ignore the $y$ dependence and thus set $N_y = 1$. The initial condition is set as a multimodal function $p_0(x) = C\{1 + 0.9\sin^2(\pi x)\}$ that is symmetric with respect to $x = 1/2$, where $C > 0$ is a normalization constant chosen so that $\int_{x\in\Omega} p_0(x)\mathrm{d}x = 1$. Using other similar $p_0$ with multiple global maxima yields results similar to those presented below. We set the parameter $q = 2$ so that the uniform distribution $p \equiv 1$ becomes an unstable equilibrium [65].

**Fig. 6** shows that the probability density $p$ develops spikes at the peaks of $p_0$, with their magnitude increasing with respect to $\delta^{-1}$. The numerical solutions are symmetric with respect to $x = 1/2$, and there exists some critical value of $\delta = \delta_c \in (0.402, 0.403)$ such that the density in the RD develops spikes, while that in the MFG does not. The two become close when $\delta = 0.403$ or larger. This suggests that the MFG solution is more robust to perturbations in a wider range of $\delta$. The formation of spikes in the RD is attributed to the positive feedback encoded in the utility, which now reads $U(x, p) = p$. This formulation implies that concentrating the mass of probability increases utility, amplifying the positive peaks in $p_0$ into spikes. The suppression of spikes formation in the MFG case, at least for the examined values of $\delta$, suggests that less myopic agents perceive the cost of switching actions to be higher.

We also compute the non-exchangeable case with the initial condition $p_0(x, y) = C\{1 + 0.9\sin^2(4\pi x)\}$, where $C$ is a normalization function such that $\int_{x\in\Omega} p_0(x, y)\mathrm{d}x = 1$ ($y \in I$). The utility is designed as follows to inherit the structure of the exchangeable case while accounting for heterogeneity among agents:

$$U(x, y, p) = \frac{yp(x)}{1 + 10A(\alpha)} \text{ with } \alpha = \int_{\Omega\times I} xp(x, y)\mathrm{d}x\mathrm{d}y\,. \tag{39}$$

The discretization of this utility is performed following (36). According to (39), a larger average of the average $\alpha$ leads to a smaller utility, similar to the tragedy of commons case considered in the previous subsection.

**Fig. 7** compares the computed probability density $p$ between the RD and MFG for the *q*-voter model. The computational results again suggest that the RD is more sensitive to the non-uniformity of the solution profile because its probability density develops peaks at the maximum points in the domain when $\delta$ is sufficiently small. Qualitatively, the same stability behavior observed in the exchangeable case is seen in the non-exchangeable case for the RD, with spikes appearing in the probability density $p$, now located near the boundary $y = 1$. More specifically, these spikes are situated at some of the grid points $P_{\cdot,N_y}$, where the utility is maximized with respect to $y$. The actions of the most aggressive agents are concentrated on these points in this computational case. The critical value of the discount rate for spike creation in the RDs in this case is $\delta_c \in (0.400, 0.401)$. The spiking phenomenon observed near the critical $\delta$ are comparable between **Figs. 6 and 7**. The computational results obtained in this subsection may contribute to better understand the behavior of solutions to the *q*-voter model with continuous or a large degree of actions.

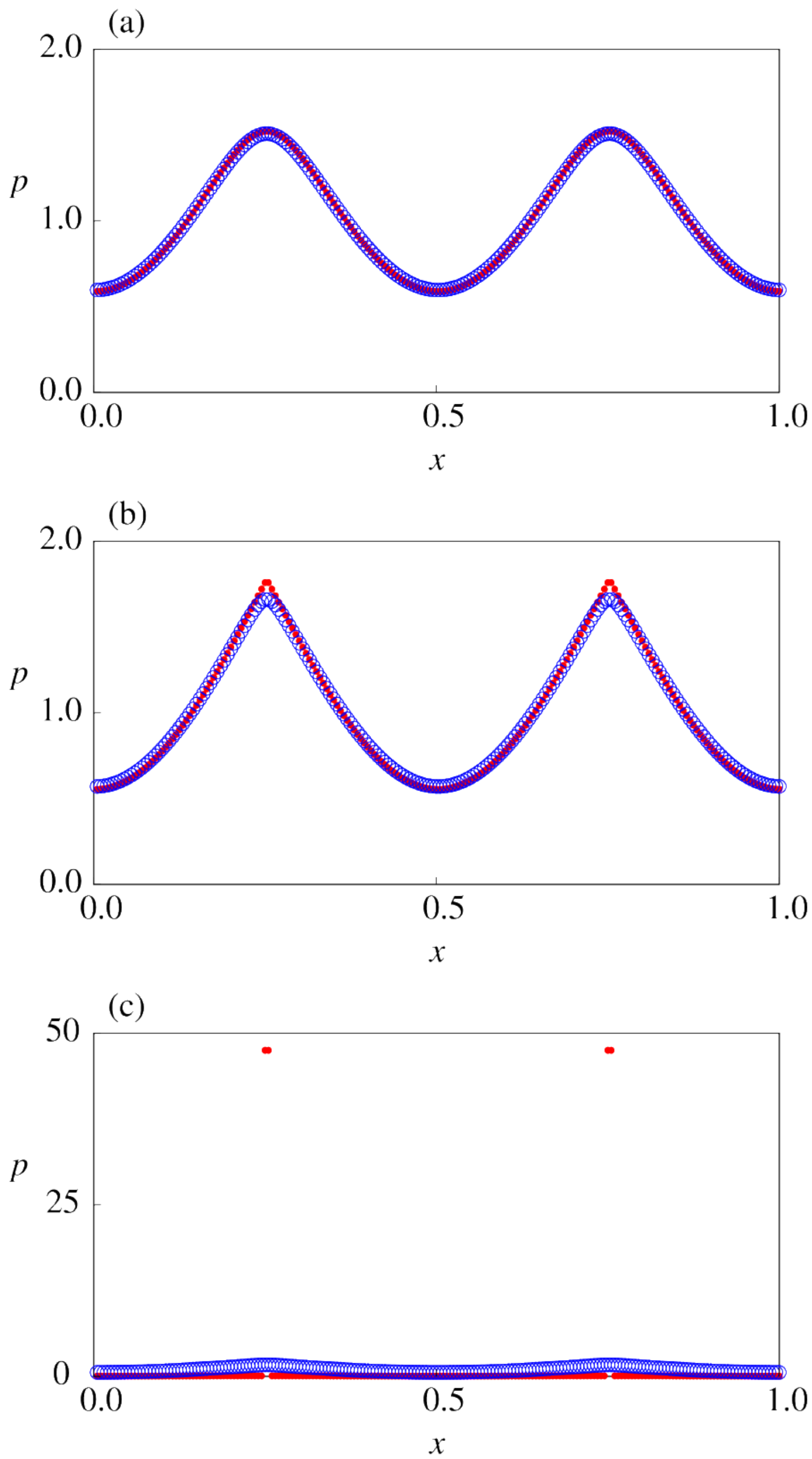

**Fig. 6** Computed probability density $p$ for the non-exchangeable $q$-voter model (Red: RD, Blue: MFG). (a) $\delta = 0.5$, (b) $\delta = 0.403$, and (c) $\delta = 0.402$.

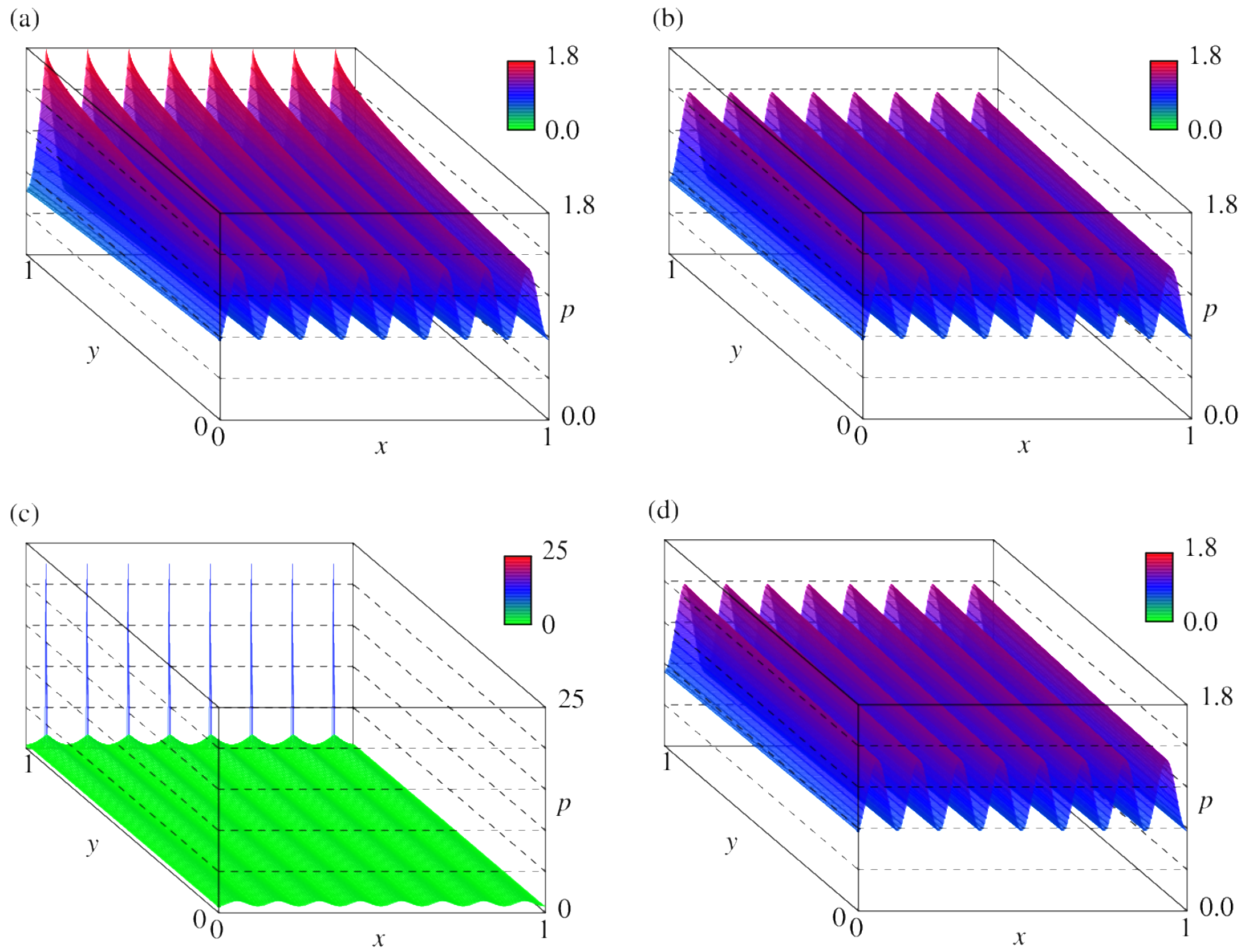

**Fig. 7** Computed probability density $p$ for the exchangeable $q$-voter model. (a) RD with $\delta = 0.401$, (b) MFG with $\delta = 0.401$, (c) RD with $\delta = 0.400$, and (d) MFG with $\delta = 0.400$.

## 5 Conclusion

In this paper, we formulated an RD and its stationary MFG version, both played by heterogeneous agents in a continuous action space. Both models were shown to be well-posed under certain conditions. Specifically, the Lipschitz continuity of the utility with respect to probability measures and the magnitude of discount rate were critical in obtaining these results. We also proposed adapting the $q$-voter model to our framework, providing a new perspective on this model in terms of a pairwise comparison protocol. The computational results for both the RD and MFG demonstrated that heterogeneity qualitatively affects the solutions. The tragedy of the commons case showed that introducing discontinuous penalization can potentially prevent excessive resource extraction. In the $q$-voter case, we found that the absence of regularity assumptions in the utility leads to singular distributions of agents' actions, where the location of the singularities depends on the agent's heterogeneity. The connection between the RD and MFG found in this paper is theoretically important because it gives a new interpretation for both of them, which would potentially advance the establishment of new mathematical theories and numerical algorithms for dynamic games. Moreover, we consider that the RD and MFG can be used for studying applied problems other than that studied in this paper, provided that the utility is suitably designed.

A limitation of this study is the absence of micro-macro linkage in theory, where the MFG should ideally be derived from a microscopic system containing a large number of agents—something that was implicitly assumed in this paper. Addressing this issue would require incorporating notions of the propagation of chaos (the mixing of actions among agents), which would necessitate a more sophisticated theoretical framework. Additionally, the non-exchangeability among agents may pose a challenge to this approach. The existence and uniqueness of solutions to the MFG when the discount rate is small remains an open question. Resolving this issue will likely require a mathematical tool beyond the fixed-point argument, which heavily depends on the discount rate. From an application perspective, a limitation is the full determination of the utility, which we plan to address in future studies focused on energy and environmental management contexts. The proposed approach could extend to MFGs with a finite horizon, provided the discount is large. This will be investigated in future research, particularly in problems involving seasonality. Moreover, applying the proposed framework to the other protocols, such as logit and projection-based ones, represents another interesting avenue for future exploration. We consider that multi-dimensional agents' dynamics have a larger degree-of-freedom than the cases we discussed in this paper, which would be an interesting future topic. Studying cases with discrete state space and/or action will also be an interesting topic.

## Appendix

### A.1 Auxiliary convergence results

We present the convergence results from the RD to the corresponding MFG for the model discussed in **Section 4.3**. We fix $P = 1$. **Table A1** shows the maximum absolute difference of the probability density $p$ between the RD and MFG. **Table A2** shows the maximum and average absolute differences between the utility $U$ of the RD and the value function $\Phi$ of the MFG. As shown in these tables, the difference between the two models decreases as $\delta$ increases even with the discontinuous utility.

**Table A1** Maximum and average absolute differences of the probability density $p$ between RD and MFG for the problem in **Section 4.3**.

| $\delta$ | Maximum difference | Average difference |
|---|---|---|
| 0.01 | $1.913 \times 10^{2}$ | $2.067 \times 10^{-1}$ |
| 0.1 | $1.407 \times 10^{2}$ | $3.030 \times 10^{-1}$ |
| 1 | $2.625 \times 10^{0}$ | $6.278 \times 10^{-2}$ |
| 10 | $2.716 \times 10^{-2}$ | $3.523 \times 10^{-3}$ |

**Table A2** Maximum and average absolute differences between the utility $U$ of RD and the value function $\Phi$ of MFG for the problem in **Section 4.3**.

| $\delta$ | Maximum difference | Average difference |
|---|---|---|
| 0.01 | $3.105 \times 10^{0}$ | $6.567 \times 10^{-1}$ |
| 0.1 | $2.710 \times 10^{0}$ | $5.027 \times 10^{-1}$ |
| 1 | $1.835 \times 10^{0}$ | $2.492 \times 10^{-1}$ |
| 10 | $6.076 \times 10^{-1}$ | $5.743 \times 10^{-2}$ |

### A.2 Proofs of Propositions

***Proof of Proposition 1***

The proof is a version of that of Theorem 12 in Mendoza-Palacios and Hernández-Lerma [53]. Therefore, we provide only a sketch of the proof. We begin by considering the following auxiliary equation:

$$\frac{\mathrm{d}\mu_t(y,\mathrm{d}x)}{\mathrm{d}t} = \left(2-\|\mu_t\|_{\mathrm{TV}}^{(l)}\right)_+ \left(U(x,y,\mu_t) - \int_{z\in\Omega} U(z,y,\mu_t)\mu_t(y,\mathrm{d}z)\right)\mu_t(y,\mathrm{d}x) \quad \left(=\mathbb{G}(\mu_t)\right), \tag{40}$$

which is an RD having a modified right-hand side multiplied by the coefficient $\left(2-\|\mu_t\|_{\mathrm{TV}}^{(l)}\right)_+$. We see that $\mathbb{G}$ is strictly bounded and Lipschitz continuous with respect to $\mu \in \mathcal{M}^{(l)}$, showing that the modified equation (40) admits a unique solution $\mu \in C\left([0,T],\mathcal{M}^{(l)}\right)$ (the left-hand side of (40) is understood in the strong sense under the TV norm) (Theorem 10, Proposition 11, and the paragraph just after Proposition 11 in Mendoza-Palacios and Hernández-Lerma [53]). This solution must be nonnegative because $\mu_0 \in \mathcal{P}^{(l)}$ and the right-hand side is proportional to $\mu_t$. Moreover, we have

$$\begin{aligned}
&\frac{\mathrm{d}}{\mathrm{d}t}\int_{x\in\Omega}\mu_t(y,\mathrm{d}x)\\
&=\int_{x\in\Omega}\frac{\mathrm{d}\mu_t(y,\mathrm{d}x)}{\mathrm{d}t}\\
&=\left(2-\|\mu_t\|_{\mathrm{TV}}^{(l)}\right)_+\int_{x\in\Omega}\left(U(x,y,\mu_t)-\int_{z\in\Omega}U(z,y,\mu_t)\mu_t(y,\mathrm{d}z)\right)\mu_t(y,\mathrm{d}x)\\
&=\left(2-\|\mu_t\|_{\mathrm{TV}}^{(l)}\right)_+\left\{\int_{x\in\Omega}U(x,y,\mu_t)\mu_t(y,\mathrm{d}x)-\int_{x\in\Omega}\int_{z\in\Omega}U(z,y,\mu_t)\mu_t(y,\mathrm{d}z)\mu_t(y,\mathrm{d}x)\right\}\\
&=\left(2-\|\mu_t\|_{\mathrm{TV}}^{(l)}\right)_+\left\{\left(\int_{z\in\Omega}U(z,y,\mu_t)\mu_t(y,\mathrm{d}z)\right)-\int_{x\in\Omega}\left(\int_{z\in\Omega}U(z,y,\mu_t)\mu_t(y,\mathrm{d}z)\right)\mu_t(y,\mathrm{d}x)\right\}\\
&=\left(2-\|\mu_t\|_{\mathrm{TV}}^{(l)}\right)_+\left(1-\int_{x\in\Omega}\mu_t(y,\mathrm{d}x)\right)\int_{z\in\Omega}U(z,y,\mu_t)\mu_t(y,\mathrm{d}z)
\end{aligned}\quad, \tag{41}$$

demonstrating that the modified equation has a critical point when $\int_{x\in\Omega}\mu_t(y,\mathrm{d}x)=1$, i.e., if $\mu_t(y,\cdot)$ is a probability measure. Moreover, its non-negativity follows from the fact that the right-hand side of (41) is proportional to $\mu_t(y,\mathrm{d}x)$ multiplied by a strictly bounded coefficient. Therefore, the modified equation admits a unique solution $\mu \in C\left([0,T],\mathcal{P}^{(l)}\right)$ for each $T>0$. This solution is also the unique solution to the RD because of the Lipschitz continuity of the right-hand side of (4) for each $\mu_t \in \mathcal{M}_2^{(l)}$.

□

***Proof of Proposition 2***

Fix some $l \in \mathbb{W}$. We first show the *a priori* estimate of solutions $\Phi \in \mathbb{X}$ to the HJB equation. Because $U$ is non-negative, we have

$$\Phi(x,y) \geq \frac{1}{2\delta}\int_{\Omega}\left(\Phi(z,y)-\Phi(x,y)\right)_+^2 m(y,\mathrm{d}z) \geq 0\,. \tag{42}$$

Assume that $\Phi$ is maximized at some $(x', y') \in \Omega \times I$ (such a point exists because $\Omega \times I$ is a compact set). Then, we have

$$\begin{aligned}\Phi(x', y') &= \mathbb{H}_l(\Phi)(x', y') \\ &\le \bar{U} + \frac{1}{2\delta}\int_\Omega \left(\Phi(z, y') - \Phi(x', y')\right)_+^2 m(y, \mathrm{d}z) \\ &= \bar{U} + 0 \\ &= \bar{U}\end{aligned}. \tag{43}$$

Therefore, any solution $\Phi \in \mathbb{X}$ to the HJB equation must belong to $\mathbb{X}_{\bar{U}}$.

Next, we consider the auxiliary equation, which is a regularized version of the HJB equation:

$$\Phi(x, y) = U(x, y, l) + \frac{1}{2\delta}\int_\Omega \left(\hat{\Phi}(z, y) - \hat{\Phi}(x, y)\right)_+^2 l(y, \mathrm{d}z) \ \left(= \hat{\mathbb{H}}_l(\Phi)(x, y)\right), \tag{44}$$

where we formally write $\hat{\Phi} = \max\left\{0, \min\left\{\bar{U}, \Phi\right\}\right\}$. Note that the *a priori* bound $0 \le \Phi \le \bar{U}$ also holds true for the auxiliary equation (44).

We show that $\hat{\mathbb{H}}_l$ is a mapping from $\mathbb{X}$ to $\mathbb{X}$. For each $\Phi \in \mathbb{X}$, we have that $\hat{\mathbb{H}}_l(\Phi)$ is a continuous function on $\Omega \times I$ because of the assumptions $l \in \mathbb{W}$ and $\Phi \in \mathbb{X}$. The boundedness of $\hat{\mathbb{H}}_l(\Phi)$ in $\mathbb{X}$ is proven as follows: at each $(x, y) \in \Omega \times I$,

$$\begin{aligned}\left|\hat{\mathbb{H}}_l(\Phi)(x, y)\right| &\le \left|U(x, y, l)\right| + \frac{1}{2\delta}\int_\Omega \left(\hat{\Phi}(z, y) - \hat{\Phi}(x, y)\right)_+^2 l(y, \mathrm{d}z) \\ &\le \bar{U} + \frac{1}{2\delta}\int_\Omega \bar{U}^2 l(y, \mathrm{d}z) \\ &\le \bar{U} + \frac{\bar{U}^2}{2\delta} \\ &< +\infty\end{aligned}. \tag{45}$$

We also have that $\hat{\mathbb{H}}_l$ is strictly contractive when $\delta > 2\bar{U}$; for each $\Phi_1, \Phi_2 \in \mathbb{X}$, we have

$$
\begin{aligned}
&\left|\hat{\mathbb{H}}_l(\Phi_1)(x,y)-\hat{\mathbb{H}}_l(\Phi_2)(x,y)\right| \\
&=\frac{1}{2\delta}\left|\int_\Omega\left(\hat{\Phi}_1(z,y)-\hat{\Phi}_1(x,y)\right)_+^2 l(y,\mathrm{d}z)-\int_\Omega\left(\hat{\Phi}_2(z,y)-\hat{\Phi}_2(x,y)\right)_+^2 l(y,\mathrm{d}z)\right| \\
&\leq\frac{1}{2\delta}\int_\Omega\left|\left(\hat{\Phi}_1(z,y)-\hat{\Phi}_1(x,y)\right)_+^2-\left(\hat{\Phi}_2(z,y)-\hat{\Phi}_2(x,y)\right)_+^2\right| l(y,\mathrm{d}z) \\
&=\frac{1}{2\delta}\int_\Omega\left(\begin{array}{l}\left|\left(\hat{\Phi}_1(z,y)-\hat{\Phi}_1(x,y)\right)_+ +\left(\hat{\Phi}_2(z,y)-\hat{\Phi}_2(x,y)\right)_+\right| \\ \times\left|\left(\hat{\Phi}_1(z,y)-\hat{\Phi}_1(x,y)\right)_+ -\left(\hat{\Phi}_2(z,y)-\hat{\Phi}_2(x,y)\right)_+\right|\end{array}\right) l(y,\mathrm{d}z) \\
&\leq\frac{1}{2\delta}\int_\Omega 2\bar{U}\left|\left(\hat{\Phi}_1(z,y)-\hat{\Phi}_1(x,y)\right)_+ -\left(\hat{\Phi}_2(z,y)-\hat{\Phi}_2(x,y)\right)_+\right| l(y,\mathrm{d}z) \\
&\leq\frac{1}{2\delta}\int_\Omega 2\bar{U}\left|\hat{\Phi}_1(z,y)-\hat{\Phi}_2(z,y)+\hat{\Phi}_2(x,y)-\hat{\Phi}_1(x,y)\right| l(y,\mathrm{d}z) \\
&\leq\frac{\bar{U}}{\delta}\int_\Omega 2\left\|\hat{\Phi}_1-\hat{\Phi}_2\right\|_{\mathbb{X}} l(y,\mathrm{d}z) \\
&=\frac{2\bar{U}}{\delta}\left\|\hat{\Phi}_1-\hat{\Phi}_2\right\|_{\mathbb{X}} \\
&\leq\frac{2\bar{U}}{\delta}\left\|\Phi_1-\Phi_2\right\|_{\mathbb{X}}
\end{aligned}
\quad . \tag{46}
$$

Therefore, we obtain the estimate

$$
\left\|\hat{\mathbb{H}}_l(\Phi_1)-\hat{\mathbb{H}}_l(\Phi_2)\right\|_{\mathbb{X}}\leq\frac{2\bar{U}}{\delta}\left\|\Phi_1-\Phi_2\right\|_{\mathbb{X}}, \tag{47}
$$

and strict contraction property follows owing to $\delta>2\bar{U}$. Then, the Banach's fixed point theorem (Theorem 5.7 in Brezis [15]) combined with the *a priori* estimate shows that there exists a unique solution $\Phi'\in\mathbb{X}_{\bar{U}}$ to the auxiliary equation (44). This $\Phi'$ also solves the equation $\Phi=\mathbb{H}_l(\Phi)$; moreover, such a solution to this equation is unique. Indeed, if there exist two solutions $\Phi_1,\Phi_2\in\mathbb{X}_{\bar{U}}$ to it, then we must have

$$
\left\|\Phi_1-\Phi_2\right\|_{\mathbb{X}}\leq\left\|\mathbb{H}_l(\Phi_1)-\mathbb{H}_l(\Phi_2)\right\|_{\mathbb{X}}\leq\frac{2\bar{U}}{\delta}\left\|\Phi_1-\Phi_2\right\|_{\mathbb{X}}<\left\|\Phi_1-\Phi_2\right\|_{\mathbb{X}}, \tag{48}
$$

and hence the uniqueness $\left\|\Phi_1-\Phi_2\right\|_{\mathbb{X}}=0$.

Finally, the unique solution to the equation $\Phi=\mathbb{H}_{l_i}(\Phi)$ with $l=l_i\in\mathbb{W}$ is denoted as $\Phi_i$ ($i=1,2$). At each $(x,y)\in\Omega\times I$, we have

$$
\begin{aligned}
&\left|\Phi_1(x,y)-\Phi_2(x,y)\right| \\
&=\left|\hat{\mathbb{H}}_{l_1}(\Phi_1)(x,y)-\hat{\mathbb{H}}_{l_2}(\Phi_2)(x,y)\right| \\
&\leq \left|U(x,y,l_1)-U(x,y,l_2)\right| \\
&\quad+\frac{1}{2\delta}\left|\int_\Omega\left(\Phi_1(z,y)-\Phi_1(x,y)\right)_+^2 l_1(y,\mathrm{d}z)-\int_\Omega\left(\Phi_2(z,y)-\Phi_2(x,y)\right)_+^2 l_2(y,\mathrm{d}z)\right| \\
&\leq L_U\left\|l_1-l_2\right\|_{\mathbb{Y}}+\frac{1}{2\delta}\left(\begin{array}{l}\left|\int_\Omega\left(\Phi_1(z,y)-\Phi_1(x,y)\right)_+^2 l_1(y,\mathrm{d}z)-\int_\Omega\left(\Phi_1(z,y)-\Phi_1(x,y)\right)_+^2 l_2(y,\mathrm{d}z)\right| \\ +\left|\int_\Omega\left(\Phi_1(z,y)-\Phi_1(x,y)\right)_+^2 l_2(y,\mathrm{d}z)-\int_\Omega\left(\Phi_2(z,y)-\Phi_2(x,y)\right)_+^2 l_2(y,\mathrm{d}z)\right|\end{array}\right) \\
&\leq L_U\left\|l_1-l_2\right\|_{\mathbb{Y}}+\frac{1}{2\delta}\left(\begin{array}{l}\left|\int_\Omega\left(\Phi_1(z,y)-\Phi_1(x,y)\right)_+^2\left(l_1(y,\mathrm{d}z)-l_2(y,\mathrm{d}z)\right)\right| \\ +\int_\Omega\left|\left(\Phi_1(z,y)-\Phi_1(x,y)\right)_+^2-\left(\Phi_2(z,y)-\Phi_2(x,y)\right)_+^2\right| l_2(y,\mathrm{d}z)\end{array}\right) \\
&\leq L_U\left\|l_1-l_2\right\|_{\mathbb{Y}}+\frac{1}{2\delta}\left(\bar{U}^2\left\|l_1-l_2\right\|_{\mathbb{Y}}+4\bar{U}\left\|\Phi_1-\Phi_2\right\|_{\mathbb{X}}\right)
\end{aligned}
\quad , (49)
$$

where we performed the calculation as in (46). Consequently, we have

$$
\left\|\Phi_1-\Phi_2\right\|_{\mathbb{X}} \leq\left(L_U+\frac{\bar{U}^2}{2\delta}\right)\left\|l_1-l_2\right\|_{\mathbb{Y}}+\frac{2\bar{U}}{\delta}\left\|\Phi_1-\Phi_2\right\|_{\mathbb{X}}, \tag{50}
$$

and hence the desired result (26).

□

***Proof of Proposition 3***

The storyline of the proof is similar to that of **Proof of Proposition 2**, but with several technical differences between them.

Fix some $\theta \in \mathbb{X}_{\bar{U}}$. We consider the following auxiliary equation, which is a regularized version of the FP equation:

$$
m(y,\mathrm{d}x)=\frac{\delta}{\delta-\left(\theta(x,y)-\max\left\{0,\int_{z\in\Omega}\theta(z,y)m(y,\mathrm{d}z)\right\}\right)}\mu_0(y,\mathrm{d}x)\ \left(=\hat{\mathbb{F}}_\theta(m)(x,y)\right). \tag{51}
$$

The coefficient multiplied by $\mu_0(y,\mathrm{d}x)$ is a bounded and continuous function on $\Omega\times I$ because its denominator is continuous and strictly positive:

$$
\begin{aligned}
\delta-\left(\theta(x,y)-\max\left\{0,\int_{z\in\Omega}\theta(z,y)m(y,\mathrm{d}z)\right\}\right)&=\delta-\theta(x,y)+\max\left\{0,\int_{z\in\Omega}\theta(z,y)m(y,\mathrm{d}z)\right\} \\
&\geq\delta-\theta(x,y) \\
&\geq\delta-\bar{U} \\
&>0
\end{aligned}
\quad . (52)
$$

Fix some $m\in\mathbb{Y}$. By (52), we have the following boundedness property (we use $\left\|\mu_0\right\|_{\mathbb{Y}}=1$):

$$
\left\|\mathbb{F}_\Phi(m)\right\|_{\mathbb{Y}} \leq\left|\frac{\delta}{\delta-\left(\theta(x,y)-\max\left\{0,\int_{z\in\Omega}\theta(z,y)m(y,\mathrm{d}z)\right\}\right)}\right|\left\|\mu_0\right\|_{\mathbb{Y}} \leq\frac{\delta}{\delta-\bar{U}}. \tag{53}
$$

Fix $m_1, m_2 \in \mathbb{Y}$. We have the following strict contraction property when $\delta > 0$ is sufficiently large:

$$
\begin{aligned}
&\left\| \mathbb{F}_{\Phi}(m_1) - \mathbb{F}_{\Phi}(m_2) \right\|_{\mathbb{Y}} \\
&= \left\| \left\{ \begin{array}{l} \dfrac{\delta}{\delta - \left(\theta(x,y) - \max\left\{0, \int_{z\in\Omega} \theta(z,y) m_1(y,\mathrm{d}z)\right\}\right)} \\ -\dfrac{\delta}{\delta - \left(\theta(x,y) - \max\left\{0, \int_{z\in\Omega} \theta(z,y) m_2(y,\mathrm{d}z)\right\}\right)} \end{array} \right\} \mu_0(y;\mathrm{d}x) \right\|_{\mathbb{Y}} \\
&\le \frac{\delta}{\left(\delta - \bar{U}\right)^2} \left\| \begin{array}{l} \delta - \left(\theta(x,y) - \max\left\{0, \int_{z\in\Omega} \theta(z,y) m_1(y,\mathrm{d}z)\right\}\right) \\ -\left(\delta - \left(\theta(x,y) - \max\left\{0, \int_{z\in\Omega} \theta(z,y) m_2(y,\mathrm{d}z)\right\}\right)\right) \end{array} \right\|_{\mathbb{X}} \left\| \mu_0 \right\|_{\mathbb{Y}} \quad . \\
&\le \frac{\delta}{\left(\delta - \bar{U}\right)^2} \left\| \left( \int_{z\in\Omega} \theta(z,y) m_1(y,\mathrm{d}z) \right)_+ - \left( \int_{z\in\Omega} \theta(z,y) m_2(y,\mathrm{d}z) \right)_+ \right\|_{\mathbb{X}} \left\| \mu_0 \right\|_{\mathbb{Y}} \\
&\le \frac{\delta}{\left(\delta - \bar{U}\right)^2} \left\| \int_{z\in\Omega} \theta(z,y) m_1(y,\mathrm{d}z) - \int_{z\in\Omega} \theta(z,y) m_2(y,\mathrm{d}z) \right\|_{\mathbb{X}} \\
&\le \frac{\delta \bar{U}}{\left(\delta - \bar{U}\right)^2} \left\| m_1 - m_2 \right\|_{\mathbb{Y}}
\end{aligned} \tag{54}
$$

We have $\dfrac{\delta \bar{U}}{\left(\delta - \bar{U}\right)^2} \in (0,1)$ owing to $\delta > \dfrac{3+\sqrt{5}}{2}\bar{U}$. Then, for such $\delta$, the Banach's fixed point theorem (Theorem 5.7 in Brezis [15]) shows that there exists a unique solution $m' \in \mathbb{Y}$ to the equation $m = \mathbb{F}_{\theta}(m)$. This solution is a non-negative measure, which directly follows from (51). Moreover, the solution is a probability measure because of the equality

$$
\int_{x\in\Omega} \left\{ \delta - \left( \theta(x,y) - \int_{z\in\Omega} \theta(z,y) m(y,\mathrm{d}z) \right) \right\} m(y,\mathrm{d}x) = \delta \int_{x\in\Omega} \mu_0(y,\mathrm{d}x) = \delta \,, \tag{55}
$$

leading to

$$
\delta \int_{x\in\Omega} m(y,\mathrm{d}x) - \delta + \int_{z\in\Omega} \theta(z,y) m(y,\mathrm{d}z) \int_{x\in\Omega} m(y,\mathrm{d}x) - \int_{x\in\Omega} \theta(z,y) m(y,\mathrm{d}z) = 0 \tag{56}
$$

and hence

$$
\left( \delta + \int_{z\in\Omega} \theta(z,y) m_2(y,\mathrm{d}z) \right) \left( \int_{x\in\Omega} m(y,\mathrm{d}x) - 1 \right) = 0 \,, \tag{57}
$$

yielding $\int_{x\in\Omega} m(y,\mathrm{d}x) = 1$ (recall that $m$ is non-negative, so that the left coefficient in (57) is positive). This shows $m \in \mathbb{W}$ for the unique solution to (51). Moreover, it is immediate to check that this solution satisfies the FP equation $m = \mathbb{F}_{\theta}(m)$. We show that the FP equation in this case admits a unique solution in $\mathbb{W}$. Indeed, if there are two solutions $m_1, m_2 \in \mathbb{W}$, then we have

$$
\left\| m_1 - m_2 \right\|_{\mathbb{Y}} = \left\| \mathbb{F}_{\theta}(m_1) - \mathbb{F}_{\theta}(m_2) \right\|_{\mathbb{Y}} \le \frac{\delta \bar{U}}{\left(\delta - \bar{U}\right)^2} \left\| m_1 - m_2 \right\|_{\mathbb{Y}} < \left\| m_1 - m_2 \right\|_{\mathbb{Y}} \,, \tag{58}
$$

showing the uniqueness.

Finally, the unique solution to equation $m=\mathbb{F}_{\theta_i}(m)$ with $\theta=\theta_i\in\mathbb{X}_{\bar{U}}$ is denoted as $m_i$ ($i=1,2$). We have

$$
\begin{aligned}
&\left\|m_1-m_2\right\|_{\mathbb{Y}}\\
&=\left\|\left\{\frac{\delta}{\delta-\left(\theta_1(x,y)-\int_{z\in\Omega}\theta_1(z,y)m_1(y,\mathrm{d}z)\right)}-\frac{\delta}{\delta-\left(\theta_2(x,y)-\int_{z\in\Omega}\theta_2(z,y)m_2(y,\mathrm{d}z)\right)}\right\}\mu_0(y;\mathrm{d}x)\right\|_{\mathbb{Y}}\\
&\leq\left\|\left\{\frac{\delta}{\delta-\left(\theta_1(x,y)-\int_{z\in\Omega}\theta_1(z,y)m_1(y,\mathrm{d}z)\right)}-\frac{\delta}{\delta-\left(\theta_1(x,y)-\int_{z\in\Omega}\theta_1(z,y)m_2(y,\mathrm{d}z)\right)}\right\}\mu_0(y;\mathrm{d}x)\right\|_{\mathbb{Y}}\\
&+\left\|\left\{\frac{\delta}{\delta-\left(\theta_1(x,y)-\int_{z\in\Omega}\theta_1(z,y)m_2(y,\mathrm{d}z)\right)}-\frac{\delta}{\delta-\left(\theta_2(x,y)-\int_{z\in\Omega}\theta_2(z,y)m_2(y,\mathrm{d}z)\right)}\right\}\mu_0(y;\mathrm{d}x)\right\|_{\mathbb{Y}}\\
&\leq\left\|\frac{\delta}{\delta-\left(\theta_1(x,y)-\int_{z\in\Omega}\theta_1(z,y)m_1(y,\mathrm{d}z)\right)}-\frac{\delta}{\delta-\left(\theta_1(x,y)-\int_{z\in\Omega}\theta_1(z,y)m_2(y,\mathrm{d}z)\right)}\right\|_{\mathbb{X}}\left\|\mu_0\right\|_{\mathbb{Y}}\\
&+\left\|\frac{\delta}{\delta-\left(\theta_1(x,y)-\int_{z\in\Omega}\theta_1(z,y)m_2(y,\mathrm{d}z)\right)}-\frac{\delta}{\delta-\left(\theta_2(x,y)-\int_{z\in\Omega}\theta_2(z,y)m_2(y,\mathrm{d}z)\right)}\right\|_{\mathbb{X}}\left\|\mu_0\right\|_{\mathbb{Y}}\\
&\leq\frac{\delta}{\left(\delta-\bar{U}\right)^2}\left\|\int_{z\in\Omega}\theta_1(z,y)m_1(y,\mathrm{d}z)-\int_{z\in\Omega}\theta_1(z,y)m_2(y,\mathrm{d}z)\right\|_{\mathbb{X}}\\
&+\frac{\delta}{\left(\delta-\bar{U}\right)^2}\left\|\theta_1(x,y)-\int_{z\in\Omega}\theta_1(z,y)m_2(y,\mathrm{d}z)-\left(\theta_2(x,y)-\int_{z\in\Omega}\theta_2(z,y)m_2(y,\mathrm{d}z)\right)\right\|_{\mathbb{X}}
\end{aligned}
\quad .(59)
$$

We also have

$$
\left\|\int_{z\in\Omega}\theta_1(z,y)m_1(y,\mathrm{d}z)-\int_{z\in\Omega}\theta_1(z,y)m_2(y,\mathrm{d}z)\right\|_{\mathbb{X}}\leq\bar{U}\left\|m_1-m_2\right\|_{\mathbb{Y}} \quad (60)
$$

and

$$
\begin{aligned}
&\left\|\theta_1(x,y)-\int_{z\in\Omega}\theta_1(z,y)m_2(y,\mathrm{d}z)-\left(\theta_2(x,y)-\int_{z\in\Omega}\theta_2(z,y)m_2(y,\mathrm{d}z)\right)\right\|_{\mathbb{X}}\\
&\leq\left\|\theta_1-\theta_2\right\|_{\mathbb{X}}+\left\|\int_{z\in\Omega}\theta_1(z,y)m_2(y,\mathrm{d}z)-\int_{z\in\Omega}\theta_2(z,y)m_2(y,\mathrm{d}z)\right\|_{\mathbb{X}}\\
&\leq\left\|\theta_1-\theta_2\right\|_{\mathbb{X}}+\left\|\int_{z\in\Omega}\left|\theta_1(z,y)-\theta_2(z,y)\right|m_2(y,\mathrm{d}z)\right\|_{\mathbb{X}}\\
&\leq\left\|\theta_1-\theta_2\right\|_{\mathbb{X}}+\left\|\theta_1-\theta_2\right\|_{\mathbb{X}}\left\|\int_{z\in\Omega}m_2(y,\mathrm{d}z)\right\|_{\mathbb{X}}\\
&=2\left\|\theta_1-\theta_2\right\|_{\mathbb{X}}
\end{aligned}
\quad . \quad (61)
$$

Consequently, we have

$$
\left\|m_1-m_2\right\|_{\mathbb{Y}}\leq\frac{\delta}{\left(\delta-\bar{U}\right)^2}\left(\bar{U}\left\|m_1-m_2\right\|_{\mathbb{Y}}+2\left\|\theta_1-\theta_2\right\|_{\mathbb{X}}\right) \quad (62)
$$

and hence the desired result (27).

□

***Proof of Proposition 4***

According to **Proposition 3**, for a given $\Phi \in \mathbb{X}_{\bar{U}}$, the solution to the FP equation (24) can be expressed as $m = \mathbb{V}(\Phi)$ with an injective mapping $\mathbb{V} : \mathbb{X}_{\bar{U}} \to \mathbb{W}$. Again, according to **Proposition 3**, $\mathbb{V}$ satisfies

$$\left\| \mathbb{V}(\Phi_1) - \mathbb{V}(\Phi_2) \right\|_{\mathbb{Y}} \le \left( 1 - \frac{\delta \bar{U}}{\left(\delta - \bar{U}\right)^2} \right)^{-1} \frac{2\delta}{\left(\delta - \bar{U}\right)^2} \left\| \Phi_1 - \Phi_2 \right\|_{\mathbb{X}} \tag{63}$$

for each $\Phi_1, \Phi_2 \in \mathbb{X}_{\bar{U}}$. Then, if we set the mapping $\mathbb{K}(\Phi) = \mathbb{H}_{\mathbb{V}(\Phi)}(\Phi)$ for any $\Phi \in \mathbb{X}_{\bar{U}}$, the HJB equation (25) can be rewritten as $\Phi = \mathbb{K}(\Phi)$. The boundedness of $\mathbb{K}$ follows as in **Proposition 2**. Moreover, for any $\Phi_1, \Phi_2 \in \mathbb{X}_{\bar{U}}$, at each $(x, y) \in \Omega \times I$, we have

$$\begin{aligned}
&\left| \mathbb{K}(\Phi_1)(x,y) - \mathbb{K}(\Phi_2)(x,y) \right| \\
&\le \left| U\left(x, y, \mathbb{V}(\Phi_1)\right) - U\left(x, y, \mathbb{V}(\Phi_2)\right) \right| \\
&\quad + \frac{1}{2\delta} \left| \int_\Omega \left( \Phi_1(z,y) - \Phi_1(x,y) \right)_+^2 \mathbb{V}(\Phi_1)(y, \mathrm{d}z) - \int_\Omega \left( \Phi_2(z,y) - \Phi_2(x,y) \right)_+^2 \mathbb{V}(\Phi_2)(y, \mathrm{d}z) \right| \\
&\le L_U \left\| \mathbb{V}(\Phi_1) - \mathbb{V}(\Phi_2) \right\|_{\mathbb{Y}} \\
&\quad + \frac{1}{2\delta} \left( \begin{aligned} &\left| \int_\Omega \left( \Phi_1(z,y) - \Phi_1(x,y) \right)_+^2 \mathbb{V}(\Phi_1)(y, \mathrm{d}z) - \int_\Omega \left( \Phi_1(z,y) - \Phi_1(x,y) \right)_+^2 \mathbb{V}(\Phi_2)(y, \mathrm{d}z) \right| \\ &+ \left| \int_\Omega \left( \Phi_1(z,y) - \Phi_1(x,y) \right)_+^2 \mathbb{V}(\Phi_2)(y, \mathrm{d}z) - \int_\Omega \left( \Phi_2(z,y) - \Phi_2(x,y) \right)_+^2 \mathbb{V}(\Phi_2)(y, \mathrm{d}z) \right| \end{aligned} \right) \\
&\le L_U \left\| \mathbb{V}(\Phi_1) - \mathbb{V}(\Phi_2) \right\|_{\mathbb{Y}} \\
&\quad + \frac{1}{2\delta} \left( \begin{aligned} &\left| \int_\Omega \left( \Phi_1(z,y) - \Phi_1(x,y) \right)_+^2 \left( \mathbb{V}(\Phi_1)(y, \mathrm{d}z) - \mathbb{V}(\Phi_2)(y, \mathrm{d}z) \right) \right| \\ &+ \int_\Omega \left| \left( \Phi_1(z,y) - \Phi_1(x,y) \right)_+^2 - \left( \Phi_2(z,y) - \Phi_2(x,y) \right)_+^2 \right| \mathbb{V}(\Phi_2)(y, \mathrm{d}z) \end{aligned} \right) \\
&\le L_U \left\| \mathbb{V}(\Phi_1) - \mathbb{V}(\Phi_2) \right\|_{\mathbb{Y}} + \frac{1}{2\delta} \left( \bar{U}^2 \left\| \mathbb{V}(\Phi_1) - \mathbb{V}(\Phi_2) \right\|_{\mathbb{Y}} + 4\bar{U} \left\| \Phi_1 - \Phi_2 \right\|_{\mathbb{X}} \right) \\
&= \left( L_U + \frac{\bar{U}^2}{2\delta} \right) \left\| \mathbb{V}(\Phi_1) - \mathbb{V}(\Phi_2) \right\|_{\mathbb{Y}} + \frac{2\bar{U} \left\| \Phi_1 - \Phi_2 \right\|_{\mathbb{X}}}{\delta}
\end{aligned} \tag{64}$$

Taking the supremum with respect to $(x, y) \in \Omega \times I$ in (64), and substituting (63) into it yields

$$\left\| \mathbb{K}(\Phi_1)(x,y) - \mathbb{K}(\Phi_2)(x,y) \right\|_{\mathbb{X}} \le C_\delta \left\| \Phi_1 - \Phi_2 \right\|_{\mathbb{X}} \tag{65}$$

with

$$C_\delta = \left( L_U + \frac{\bar{U}^2}{2\delta} \right) \left( 1 - \frac{\delta \bar{U}}{\left(\delta - \bar{U}\right)^2} \right)^{-1} \frac{2\delta}{\left(\delta - \bar{U}\right)^2} + \frac{2\bar{U}}{\delta}. \tag{66}$$

We can see that $C_\delta < 1$ for a sufficiently large $\delta > \frac{3+\sqrt{5}}{2}$ because $C_\delta = O\left(\delta^{-1}\right)$ for large $\delta > 0$. For such a large $\delta > 0$, the uniqueness of solutions to the HJB equation follows because of

$$\left\| \Phi_1 - \Phi_2 \right\|_{\mathbb{X}} \le \left\| \mathbb{K}(\Phi_1)(x,y) - \mathbb{K}(\Phi_2)(x,y) \right\|_{\mathbb{X}} \le C_\delta \left\| \Phi_1 - \Phi_2 \right\|_{\mathbb{X}} < \left\| \Phi_1 - \Phi_2 \right\|_{\mathbb{X}}. \tag{67}$$

Finally, the unique solution to the FP equation is then given by $m = \mathbb{V}(\Phi)$.

Finally, the fact that $m$ is absolutely continuous with respect to $\mu_0$ is the direct consequence of the fact that $m$ is given by $\mu_0$ multiplied by a strictly positive, bounded, and continuous function.

□

***Proof of Proposition 5***

Fix arbitrary $(x,y)\in\Omega\times I$. We have $\Phi\in\mathbb{X}_{\bar{U}}$ and hence it is uniformly bounded for (large) $\delta>0$. In addition, $m\in\mathbb{W}$ for any sufficiently large $\delta>0$. This implies that, in $\mathbb{W}$ we have

$$\begin{aligned}\lim_{\delta\to+\infty} m(y,\mathrm{d}x) &= \lim_{\delta\to+\infty}\frac{\delta}{\delta-\left(\Phi(x,y)-\int_{z\in\Omega}\Phi(z,y)m(y,\mathrm{d}z)\right)}\mu_0(y,\mathrm{d}x)\\ &= \lim_{\delta\to+\infty}\frac{1}{1-\delta^{-1}\left(\Phi(x,y)-\int_{z\in\Omega}\Phi(z,y)m(y,\mathrm{d}z)\right)}\mu_0(y,\mathrm{d}x).\\ &= \mu_0(y,\mathrm{d}x)\end{aligned}\tag{68}$$

Note that the quantity $\delta^{-1}\left(\Phi(x,y)-\int_{z\in\Omega}\Phi(z,y)m(y,\mathrm{d}z)\right)$ goes to 0 irrespective to the continuity of the limit function $\lim_{\delta\to+\infty}\Phi$ because of the uniform boundedness of $\Phi$ for each large $\delta>0$. With (68) in mind, by (3), in $\mathbb{X}$ we further have

$$\begin{aligned}\lim_{\delta\to+\infty}\Phi(x,y) &= \lim_{\delta\to+\infty}\left(U(x,y,m)+\frac{1}{2\delta}\int_{\Omega}\left(\Phi(z,y)-\Phi(x,y)\right)_+^2 m(y,\mathrm{d}z)\right)\\ &= \lim_{\delta\to+\infty}U(x,y,m)+0\\ &= U(x,y,\mu_0)\end{aligned}\tag{69}$$

because of

$$\lim_{\delta\to+\infty}\left|U(x,y,m)-U(x,y,\mu_0)\right|\leq\lim_{\delta\to+\infty}\left\|m-\mu_0\right\|_{\mathbb{Y}}=0\tag{70}$$

and the following inequality owing to $\Phi\in\mathbb{X}_{\bar{U}}$ and $m\in\mathbb{W}$ for each large $\delta>0$:

$$0\leq\lim_{\delta\to+\infty}\frac{1}{2\delta}\int_{\Omega}\left(\Phi(z,y)-\Phi(x,y)\right)_+^2 m(y,\mathrm{d}z)\leq\lim_{\delta\to+\infty}\frac{\bar{U}^2}{2\delta}=0.\tag{71}$$

□

***Proof of Proposition 6***

We first prove (31) through (33). Assume that (31) through (33) holds true at some $n\in\mathbb{N}\cup\{0\}$ (they are satisfied at $n=0$). Then, we have the non-negativity

$$
\begin{aligned}
m_{i,j}^{(n+1)} &= m_{i,j}^{(n)} + \Delta t \mathbb{F}_{i,j}\left(m^{(n)}, \Phi^{(n)}\right) \\
&= m_{i,j}^{(n)} + \Delta t \left\{ \delta \mu_{0,i,j} - \left( \delta - \left( \Phi_{i,j}^{(n)} - \sum_{k=1}^{N_x} \Phi_{k,j}^{(n)} m_{k,j}^{(n)} \right) \right) m_{i,j}^{(n)} \right\} \\
&= m_{i,j}^{(n)} + \Delta t \left\{ \delta \mu_{0,i,j} - \delta m_{i,j}^{(n)} + \Phi_{i,j}^{(n)} m_{i,j}^{(n)} - \left( \sum_{k=1}^{N_x} \Phi_{k,j}^{(n)} m_{k,j}^{(n)} \right) m_{i,j}^{(n)} \right\} \\
&= \left( 1 - \delta \Delta t + \Phi_{i,j}^{(n)} \Delta t - \left( \sum_{k=1}^{N_x} \Phi_{k,j}^{(n)} m_{k,j}^{(n)} \right) \Delta t \right) m_{i,j}^{(n)} + \Delta t \delta \mu_{0,i,j} \\
&\geq \left( 1 - \delta \Delta t + 0 \Delta t - \left( \sum_{k=1}^{N_x} K m_{k,j}^{(n)} \right) \Delta t \right) m_{i,j}^{(n)} + \Delta t \delta \mu_{0,i,j} \\
&\geq \left( 1 - \left( \delta + K \right) \Delta t \right) m_{i,j}^{(n)} + \Delta t \delta \mu_{0,i,j} \\
&\geq 0
\end{aligned}
\tag{72}
$$

and the conservation property

$$
\begin{aligned}
\sum_{i=1}^{N_x} m_{i,j}^{(n+1)} &= \sum_{i=1}^{N_x} m_{i,j}^{(n)} + \Delta t \sum_{i=1}^{N_x} \mathbb{F}_{i,j}\left(m^{(n)}, \Phi^{(n)}\right) \\
&= \sum_{i=1}^{N_x} m_{i,j}^{(n)} + \Delta t \sum_{i=1}^{N_x} \left\{ \delta \mu_{0,i,j} - \left( \delta - \left( \Phi_{i,j}^{(n)} - \sum_{k=1}^{N_x} \Phi_{k,j}^{(n)} m_{k,j}^{(n)} \right) \right) m_{i,j}^{(n)} \right\} \\
&= \sum_{i=1}^{N_x} m_{i,j}^{(n)} + \Delta t \sum_{i=1}^{N_x} \left\{ \delta \mu_{0,i,j} - \delta m_{i,j}^{(n)} + \Phi_{i,j}^{(n)} m_{i,j}^{(n)} - \left( \sum_{k=1}^{N_x} \Phi_{k,j}^{(n)} m_{k,j}^{(n)} \right) m_{i,j}^{(n)} \right\} . \\
&= \sum_{i=1}^{N_x} m_{i,j}^{(n)} + \Delta t \left( \delta + \sum_{k=1}^{N_x} \Phi_{k,j}^{(n)} m_{k,j}^{(n)} \right) \left( 1 - \sum_{i=1}^{N_x} m_{i,j}^{(n)} \right) \\
&= \sum_{i=1}^{N_x} m_{i,j}^{(n)} \\
&= 1
\end{aligned}
\tag{73}
$$

It remains to show the boundedness at $n+1$. By $\Delta t < \left( \delta + K \right)^{-1}$, the left inequality of (32) follows from

$$
\Phi_{i,j}^{(n+1)} = \left( 1 - \delta \Delta t \right) \Phi_{i,j}^{(n)} + \delta \Delta t U_{i,j}^{(n)} + \frac{1}{2\delta} \sum_{k=1}^{N_x} \left( \Phi_{k,j}^{(n)} - \Phi_{i,j}^{(n)} \right)_+^2 m_{k,j}^{(n)} \delta \Delta t \geq 0 . \tag{74}
$$

To show the right inequality of (32), we first consider the inequality

$$
0 = -u + \bar{U} + \frac{u^2}{2\delta} \tag{75}
$$

whose largest solution is $K$. By $\Delta t < \left( \delta + K \right)^{-1}$ and the definition of $K$, we have

$$
\begin{aligned}
\Phi_{i,j}^{(n+1)} &= \left(1-\delta\Delta t\right)\Phi_{i,j}^{(n)} + \delta\Delta t U_{i,j}^{(n)} + \frac{\delta\Delta t}{2\delta}\sum_{k=1}^{N_x}\left(\Phi_{k,j}^{(n)}-\Phi_{i,j}^{(n)}\right)_+^2 m_{k,j}^{(n)} \\
&\le \left(1-\delta\Delta t\right)K + \delta\Delta t\bar{U} + \frac{\delta\Delta t}{2\delta}\sum_{k=1}^{N_x}K^2 m_{k,j}^{(n)} \\
&= \left(1-\delta\Delta t\right)K + \delta\Delta t\bar{U} + \frac{\delta\Delta t}{2\delta}K^2 \\
&= K + \delta\Delta t\left(-K+\bar{U}+\frac{K^2}{2\delta}\right) \\
&= K
\end{aligned}
\quad . \tag{76}
$$

An induction argument combining (72) through (76).

Second, (32). The left inequality follows owing to

$$
\Phi_{i,j}^{(\infty)} = U_{i,j}^{(\infty)} + \frac{1}{2\delta}\sum_{k=1}^{N_x}\left(\Phi_{k,j}^{(\infty)}-\Phi_{i,j}^{(\infty)}\right)_+^2 m_{k,j}^{(\infty)} \ge 0 , \tag{77}
$$

where we set $\lim_{n\to+\infty} m_{i,j}^{(n)} = m_{i,j}^{(\infty)}$ that should satisfy (31) and (33) by the induction argument. The right inequality of (32) is proven as follows. Set $\bar{\Phi} = \sup_{\substack{1\le i\le N_x \\ 1\le j\le N_y}} \Phi_{i,j}^{(\infty)}$ and its maximizer $(i,j)=(\bar{i},\bar{j})$. Then, we have

$$
\begin{aligned}
\Phi_{\bar{i},\bar{j}}^{(\infty)} &= \bar{\Phi} \\
&= U_{\bar{i},\bar{j}}^{(\infty)} + \frac{1}{2\delta}\sum_{k=1}^{N_x}\left(\Phi_{k,\bar{j}}^{(\infty)}-\bar{\Phi}\right)_+^2 m_{k,\bar{j}}^{(\infty)} \\
&\le \bar{U} + \frac{1}{2\delta}\sum_{k=1}^{N_x}\left(\bar{\Phi}-\bar{\Phi}\right)_+^2 m_{k,\bar{j}}^{(\infty)} \\
&= \bar{U}
\end{aligned}
\quad . \tag{78}
$$

□

**Proof of Proposition 7**

The proof is a discrete-state analogue of those for **Propositions 2** through **4**. Therefore, we only present its sketch here. Fix arbitrary $N_x, N_y \in \mathbb{N}$, and set $N = N_x N_y$. We use the bold notations $\mathbf{\Phi} = \left\{\Phi_{i,j}\right\}_{\substack{1\le i\le N_x \\ 1\le j\le N_y}}$ and $\mathbf{m} = \left\{m_{i,j}\right\}_{\substack{1\le i\le N_x \\ 1\le j\le N_y}}$. Set the function spaces $\mathbb{X}_N = \left\{\mathbf{\Phi} \middle| \Phi_{i,j} \in \mathbb{R}^N\right\}$, $\mathbb{X}_{N,\bar{U}} = \left\{\mathbf{\Phi} \middle| \mathbf{\Phi} \in \mathbb{X}_N, 0 \le \Phi_{i,j} \le \bar{U}\right\}$, and $\mathbb{W}_N = \left\{\mathbf{m} \middle| m_{i,j} \ge 0, \sum_{i=1}^{N_x} m_{i,\cdot} = 1\right\}$. In the rest of this proof, the equation $\mathbb{H}_{i,j}(\cdot,\cdot) = 0$ ($1 \le i \le N_x$, $1 \le j \le N_y$) is abbreviated as $\mathbb{H}(\mathbf{m},\mathbf{\Phi}) = 0$. The same applies to $\mathbb{F}_{i,j}(\cdot,\cdot) = 0$.

For each $\mathbf{m} \in \mathbb{W}_N$, we have that the mapping $\mathbb{H}(\mathbf{m},\mathbf{\Phi})$ is bounded and strictly contractive for a sufficiently large $\delta > 0$, and hence $\mathbb{H}(\mathbf{m},\mathbf{\Phi}) = 0$ admits a unique solution in $\mathbb{X}_{N,\bar{U}}$ (more rigorously, we first use a regularized version as in the **Proof of Proposition 2**, and then we deal with the equation without it). Similarly, for each $\mathbf{\Phi} \in \mathbb{X}_{N,\bar{U}}$, we have that the mapping $\mathbb{F}(\mathbf{m},\mathbf{\Phi})$ is bounded and strictly contractive for a sufficiently large $\delta > 0$, and hence $\mathbb{F}(\mathbf{m},\mathbf{\Phi}) = 0$ admits a unique solution in $\mathbb{X}_{N,\bar{U}}$ (again, more rigorously we use a regularized version as in the **Proof of Proposition 3**). For any $\mathbf{u} \in \mathbb{X}_N$, we set the two discrete norms $\|\mathbf{u}\|_1 = \max_{1 \le j \le N_y} \sum_{i=1}^{N_x} |u_{i,j}|$ and $\|\mathbf{u}\|_\infty = \max_{\substack{1 \le i \le N_x \\ 1 \le j \le N_y}} |u_{i,j}|$. These norms are discrete counterparts of $\|\cdot\|_{\mathbb{Y}}$ and $\|\cdot\|_{\mathbb{X}}$ for the continuous case, respectively.

With these results in mind, for a sufficiently large $\delta > 0$, we obtain the following estimate:

$$\|\mathbf{m}_1 - \mathbf{m}_2\|_1 \le O(\delta)\|\mathbf{\theta}_1 - \mathbf{\theta}_2\|_\infty \text{ for each } \mathbf{\theta}_1, \mathbf{\theta}_2 \in \mathbb{X}_{N,\bar{U}}, \tag{79}$$

where $\mathbf{m}_i \in \mathbb{W}_N$ is the unique solution to $\mathbb{F}(\mathbf{m},\mathbf{\theta}_i) = 0$ ($i = 1,2$). We also have the estimate

$$\|\mathbf{\Phi}_1 - \mathbf{\Phi}_2\|_\infty \le O(1)\|\mathbf{l}_1 - \mathbf{l}_2\|_1 \text{ for each } \mathbf{l}_1, \mathbf{l}_2 \in \mathbb{W}_N, \tag{80}$$

where $\mathbf{m}_i \in \mathbb{W}_N$ is the unique solution to $\mathbb{H}(\mathbf{l}_i,\mathbf{\Phi}) = 0$ ($i = 1,2$). For each $\mathbf{\theta} \in \mathbb{X}_{N,\bar{U}}$, we can understand the equation $\mathbb{F}(\mathbf{m},\mathbf{\theta}) = 0$ as the mapping $\mathbf{m} = \mathbb{F}_N(\mathbf{\theta})$. Then, with (79) and (80) in mind, the mapping $\mathbb{H}(\mathbb{F}_N(\mathbf{\Phi}),\mathbf{\Phi})$ is bounded and strictly contractive for a sufficiently large $\delta > 0$. This shows that the system (29) and (30) admits a unique solution $(\mathbf{m},\Phi) \in \mathbb{W}_N \times \mathbb{X}_{N,\bar{U}}$, which is the desired result.

□